\documentclass[a4paper]{article}

\usepackage[english]{babel}
\usepackage[T1]{fontenc}
\usepackage{subcaption}
\usepackage{authblk}

\usepackage[a4paper,top=3cm,bottom=2cm,left=3cm,right=3cm,marginparwidth=1.75cm]{geometry}

\usepackage{amsmath,amsthm,amssymb,scrextend}
\usepackage{graphicx}
\usepackage[colorinlistoftodos]{todonotes}
\usepackage[colorlinks=true, allcolors=blue]{hyperref}
\usepackage{tikz}
\usepackage[utf8]{inputenc}
\usepackage{pgfplots} 
\usepackage{pgfgantt}
\usepackage{pdflscape}
\pgfplotsset{compat=newest} 
\pgfplotsset{plot coordinates/math parser=false} 
\newlength\fwidth
\newlength\fheight

\newcommand{\Q}{\ensuremath{\mathsf{Q}}}
\newcommand{\T}{\ensuremath{\mathsf{T}}}

\newcommand{\Real}{\mathbb{R}}

\newcommand{\dt}{\Delta t}

\newcommand{\y}{\mathbf y}
\newcommand{\f}{\mathbf f}
\newcommand{\q}{\mathbf q}

\newcommand{\R}{\mathbf R}

\theoremstyle{definition}
\newtheorem{remark}{Remark}
\newcommand{\Ra}[1]{\color{black}{#1}}
\newcommand{\Rb}[1]{\color{black}{#1}}
\newcommand{\Rc}[1]{\color{black}{#1}}
 \title{Hermite-based, One-step, Variational and Galerkin Time Integrators for Mechanical Systems}
\author{Harsh Sharma \thanks{Corresponding Author, Email Address for Correspondence: hasharma@ucsd.edu} \thanks{Postdoctoral Scholar, Department of Mechanical and Aerospace Engineering, University of California San Diego}, Mayuresh Patil \thanks{Professor of Practice, Department of Aerospace Engineering, Georgia Tech}, and Craig Woolsey \thanks{Professor, Kevin T. Crofton Department of Aerospace and Ocean Engineering, Virginia Tech}}
\begin{document}
\maketitle

\begin{abstract}
In this paper, we present two Hermite polynomial based approaches to derive one-step numerical integrators for mechanical systems. These methods are based on discretizing the configuration using Hermite polynomials which leads to numerical trajectories continuous in both configuration and velocity. First, we incorporate Hermite polynomials for time-discretization and derive one-step variational methods by discretizing the Lagrange-d'Alembert principle over a single time step. Second, we present the Galerkin approach to derive one-step numerical integrators by setting the weighted average of the residual of the equations of motion over a time step to zero. 
%
\par
We consider three numerical examples to understand the numerical performance of the one-step variational and Galerkin methods. We first study a particle in a double-well potential and compare the variational approach results with the corresponding results for the Galerkin approach. We then study the Duffing oscillator to understand the numerical behavior in presence of dissipative forces. Finally, we apply the proposed methods to a nonlinear aeroelastic system with two degrees of freedom. Both variational and Galerkin one-step methods capture conservative and nonconservative dynamics accurately with excellent energy behavior. The one-step Galerkin methods exhibit better trajectory and energy performance than the one-step variational methods and the variational integrators.
\end{abstract}
\section{Introduction}
\label{s:1}
The field of numerical integration of  ordinary differential equations (ODEs) reached scientific maturity towards the end of the 20th century when general purpose codes based on Runge-Kutta or linear multistep methods became widely available. In fact these off-the-shelf numerical integration packages have become an essential tool in studying complex dynamical behavior of nonlinear systems. Over the past decades, with growing interest in long-time simulation of dynamical systems from fields such as astronomy or molecular dynamics, the geometric properties of the governing differential equations became crucial for numerical simulation. This need for special numerical methods for certain classes of problems led to the advent of the field of geometric numerical integration (GNI) methods. Unlike the traditional methods that are designed to give computationally inexpensive solutions with small overall error, GNI methods are designed to preserve the underlying geometric structure of the governing differential equations. Motivation for applying these structure-preserving numerical methods to engineering problems ranges from preservation of qualitative features, such as invariants of motion or configuration space structure, to long-time stability - accurate numerical simulation for exponentially long times.  The essential aspects of GNI methods for ODEs are presented in \cite{hairer2006geometric,leimkuhler_reich_2005}.
\par 
Many mechanical systems of interest to engineers possess physically meaningful invariants such as momentum, energy, vorticity or rotational symmetry. Among such mechanical systems, Lagrangian/Hamiltonian systems form the most important class of ODEs in the context of structure-preserving numerical integration. One of the key features of Lagrangian/Hamiltonian mechanical systems is that the governing equations can be derived from variational principles. {\Rb{Wendlandt and Marsden \cite{wendlandt1997mechanical} derived mechanical integrators for conservative dynamical systems by discretizing Hamilton's principle, instead of discretizing the governing equations. These numerical integrators are both symplectic and momentum-preserving. These numerical integrators, because of the variational nature of their derivation, are known as variational integrators. Building on the ideas of Maeda, Moser and Veselov \cite{maeda1981extension,maeda1982lagrangian,moser1991discrete, veselov1991integrable,veselov1988integrable}, Marsden and his colleagues \cite{marsden2001discrete} developed the discrete mechanics framework and extended these methods to nonconservative \cite{kane2000variational} and/or constrained mechanical systems \cite{leyendecker2008variational}. Leok and Shingel \cite{leok2012prolongation} used Hermite polynomials in the discrete mechanics setting to derive variational integrators that approximate time derivatives of the trajectory with accuracy. Kane et al. \cite{kane1999symplectic} developed symplectic-energy-momentum integrators for conservative Lagrangian systems by imposing an additional discrete energy equation to compute the adaptive time step.}} {\Rc{Sharma et al. \cite{sharma2018energy} extended these energy-preserving, adaptive time step variational integrators to Lagrangian systems with external forcing.}}  A recent review paper by Sharma et al. \cite{sharma2020review} provides an overview of variational integrators with focus on their applications in different engineering fields.
\par
Alternatively, the structural dynamics community has used time finite elements to derive numerical integrators from the variational principle since the late 1960s. Argryis and Scharpf \cite{argyris1969finite} used cubic Hermite shape functions to formulate time finite elements and discretize Hamilton's principle to obtain numerical integrators for initial value problems. {\Rb{Baruch and Riff \cite{baruch1982hamilton} presented six different formulations of time finite elements based on Hamilton's law of varying action, and later on, Riff and Baruch \cite{riff1984time} implemented numerical integrators based on these formulations. Recently, Mergel et al. \cite{mergel2017c1} developed $\mathcal{C}^1-$continuous time integrators based on Hamilton's law of varying action for conservative systems.}} Unlike the discrete mechanics approach used for deriving variational integrators, the time finite elements approach only considers the variational principle over one time step and uses Hermite polynomials for discretizing the configuration which leads to $\mathcal{C}^1$--continuous trajectories. Apart from the accuracy and computational stability, discretization uses Hermite polynomials also leads to numerical methods that are amenable to feedback control implementation and analysis since they naturally yield configuration states and velocities. 
\par
{\Rb{Most of these time finite elements approaches based on Hamilton's principle or Hamilton's law of varying action fit into the framework of continuous Galerkin  or discontinuous Galerkin method. Hulme \cite{hulme1972one} employed the so-called continuous Galerkin method to derive one-step methods for numerically solving systems of nonlinear first-order ODEs. Betsch and Steinmann \cite{betsch2000conservation, betsch2001conservation}  applied Petrov-Galerkin based time finite elements to the Hamiltonian formulation of the $N$-Body problem and semidiscrete nonlinear elastodynamics.}}
\par 
The purpose of this paper is twofold. First, we derive Hermite polynomial based $\mathcal{C}^1$--continuous numerical integrators for mechanical systems from two different approaches. For the variational methods, we discretize the variational principle over a single fixed time step and then use the discrete mechanics framework to develop one-step variational methods. We also use Galerkin's method of weighted residuals of the governing equations over a single fixed time step to develop one-step Galerkin methods. Second, we study the numerical performance of the developed methods for three different classes of mechanical systems. We  also investigate the linear stability and the symplecticity of both one-step methods.
\par
The remainder of the paper is organized as follows. In Section \ref{s:2} we review the basic concepts from variational integrators and Galerkin methods used later in this work. In Section \ref{s:3} we present Hermite polynomial based one-step variational and Galerkin methods.  First, we introduce cubic Hermite polynomials used for time-discretization in this work. Then we utilize the discrete mechanics approach to derive the one-step variational methods for mechanical systems with external forcing. Finally, we use the Galerkin approach to derive one-step Galerkin methods. In Section \ref{s:4} we investigate the linear stability and symplectic nature of the proposed one-step methods. In Section \ref{s:5} we study three numerical examples to understand the numerical performance of the proposed one-step methods. Finally, in Section \ref{s:6} we provide concluding remarks and suggest future directions for this work.
\section{Background}
\label{s:2}
The numerical integration of mechanical systems can be approached in different ways. Traditional methods apply time-discretization directly to the governing equations of motion to obtain time-integration algorithms. Unfortunately, this approach does not account for the qualitative properties of the dynamical system. In this section, we review the basic concepts from variational mechanics and Galerkin methods used in the development of our one-step variational and Galerkin methods.
\par
Consider a time-invariant Lagrangian mechanical system with a finite-dimensional, smooth configuration manifold $\Q$, state space $\T\Q$, and Lagrangian $L: \T\Q \to \mathbb{R}$.
For such an autonomous Lagrangian system with time-independent external forcing $\textbf{f}_L (\q(t),\dot{\q}(t))$, the Lagrange-d'Alembert principle characterizes trajectories $\q(t)$ as those satisfying
\begin{equation}
\delta \ \int^{t_{\rm f}}_{t_{\rm i}} L(\q(t),\dot{\q}(t)) \ \text{d}t + \int ^{t_{\rm f}}_{t_{\rm i}} \textbf{f}_L (\q(t),\dot{\q}(t)) \cdot \delta \q \ \text{d} t  =0,
\label{eq:lda}
\end{equation}
where the first term considers the variation of the action integral and the second term accounts for the virtual work done by the external forces when the path $\q(t)$ is varied by $\delta \q(t)$. Using integration by parts and setting the variations at the endpoints equal to zero gives the forced Euler-Lagrange equations 
\begin{equation}
 \textbf{M}_{\rm eq}(\q(t),\dot{\q}(t),\ddot{\q}(t))=\frac{\partial L(\q(t),\dot{\q}(t))}{\partial \q}   - \frac{\text{d}}{\text{d}t} \left( \frac{\partial L(\q(t),\dot{\q}(t))}{\partial \dot{\q}} \right) + \textbf{f}_L(\q(t),\dot{\q}(t)) =\textbf{0},
 \label{eq:fg}
\end{equation} 
where $\textbf{M}_{\rm eq}$ denotes equations of motion written in a residual form. For a general mechanical system with a separable Lagrangian of the form $L(\q(t),\dot{\q}(t))=\frac{1}{2}\dot{\q}(t)^\top \mathbb{M}(\q)\dot{\q}(t)-U(\q)$, the equations of motion are given by
\begin{equation}
\textbf{M}_{\rm eq}(\q(t),\dot{\q}(t),\ddot{\q}(t))=\mathbb{M}(\q)\frac{\text{d}^2\q}{\text{d}t^2} + \frac{\partial U(\q)}{\partial \q} - \textbf{f}_L(\q(t),\dot{\q}(t)) = \textbf{0},   
\end{equation}
where $\mathbb{M}(\q)$ is the mass matrix and $U(\q)$ is the potential energy of the mechanical system.
\subsection{Discrete Variational Mechanics and Variational Integrators}
Variational integrators are time-integration methods that utilize concepts from discrete mechanics,
a discrete analogue to continuous-time variational mechanics. Although these methods were
originally developed for conservative dynamical systems, due to the variational nature of their
construction, these methods can be easily extended to nonconservative mechanical systems by
discretizing the corresponding Lagrange-d’Alembert principle. The basic idea is to first construct
discrete-time approximations of both the action integral and virtual work terms in \eqref{eq:lda} and then apply a discrete variational principle to obtain discrete-time trajectories of the mechanical system.
\par
Consider a discrete Lagrangian system with configuration manifold $\Q$ and discrete state space $\Q \times \Q$. For a fixed time step $\dt=\frac{t_{\rm f}-t_{\rm i}}{N}$, the discrete trajectory $\{\ \q_k  \  \}_{k=0}^N\ $ is defined by the configuration of the system at the sequence of times $\{\ t_k= t_{\rm i} + k\dt~|~k=0,...,N\}\ $. 
We introduce the discrete Lagrangian function $L_{\rm d}(\q_k,\q_{k+1})$, an approximation of the action integral along the curve from $\q_k$ to $\q_{k+1}$, which approximates the integral of the Lagrangian in the following sense
\begin{equation}
L_{\rm d}(\q_k,\q_{k+1}) \approx \int^{t_{k+1}}_{t_k} L(\q(t),\dot{\q}(t)) \  \text{d}t.
\label{eq:dl}
\end{equation}
To discretize external forcing, we define two discrete forces $\textbf{f}_{\rm d}^{\pm} :  \Q \times \Q \to \T^*\Q$ which approximate the continuous-time force integral that appears in~\eqref{eq:lda} over one time step in the following sense
\begin{equation}
 \textbf{f}_{\rm d}^+(\q_k,\q_{k+1})\cdot\delta \q_{k+1} +  \textbf{f}_{\rm d}^-(\q_k,\q_{k+1})\cdot\delta \q_{k}  \approx \int ^{t_{k+1}}_{t_k} \textbf{f}_L (\q(t), \dot{\q}(t)) \cdot \delta \q \ \text{d}t.
\end{equation}
The discrete Lagrange-d'Alembert principle seeks $\{\ \q_k  \  \}_{k=0}^N\ $ that satisfy 
\begin{equation}
\delta \sum_{k=0}^{N-1} L_{\rm d}(\q_k,\q_{k+1}) + \sum_{k=0}^{N-1}[\textbf{f}_{\rm d}^+(\q_k,\q_{k+1})\cdot\delta \q_{k+1} +  \textbf{f}_{\rm d}^-(\q_k,\q_{k+1})\cdot\delta \q_{k}] = 0,
\end{equation}
which yields the following forced discrete Euler-Lagrange equations 
\begin{equation}
\frac{\partial L_{\rm d}(\q_{k-1},\q_{k})}{\partial \q_k} + \frac{\partial L_{\rm d}(\q_k,\q_{k+1})}{\partial \q_{k}} + \textbf{f}_{\rm d}^+(\q_{k-1},\q_k) + \textbf{f}_{\rm d}^-(\q_{k},\q_{k+1}) = \textbf{0}   \ \ \ \ k=1,..., N-1.
\end{equation}
 These discrete equations can be recast in a standard time-marching form as follows
\begin{align}
\label{eqn:vif_imp}
-\frac{\partial L_{\rm d}(\q_k,\q_{k+1})}{\partial \q_k} - \textbf{f}_{\rm d}^-(\q_k,\q_{k+1}) &=\textbf{p}_k, \\
\label{eqn:vif_exp}
\textbf{p}_{k+1}&=\frac{\partial L_{\rm d}(\q_{k},\q_{k+1})}{\partial \q_{k+1}} + \textbf{f}_{\rm d}^+(\q_k,\q_{k+1}),
\end{align} 
where $\textbf{p}_k$ is the discrete momentum corresponding to the discrete configuration $\q_k$.
\par 
Due to the variational nature of their derivation, discrete trajectories from variational integrators inherit geometric properties from their continuous-time counterpart. These algorithms are ideal for long-time simulation because of their numerical stability and excellent energy behavior over exponentially long times for conservative as well as nonconservative systems.
\subsection{Galerkin Methods for Numerical Integration}
Galerkin methods are a class of methods used for converting continuous differential equation problems to discrete problems. The basic idea behind these methods is to seek approximate solutions to the differential equation in a finite-dimensional space spanned by a set of basis functions.
\par 
{\Rb{Consider a Hilbert space $U$, and a bilinear form $a : U \times V \to \Real$ which is both bounded and V-eliptic. We consider the following abstract problem posed as a weak formulation on the Hilbert spaces $U$ and $V$, namely,
\begin{equation}
Find\quad u \in U\quad  s.t \quad \quad a(u,v)=\ell(v)  \quad \forall v \in V,
\label{eq:weak}
\end{equation}
where $\ell$ is a bounded linear functional on $V$. In general, it is very rare to find an exact solution of \eqref{eq:weak} because $U$ and $V$ are infinite-dimensional. A natural approach to construct an approximate solution is to solve a finite-dimensional analogue of \eqref{eq:weak}. Let $U_N \subseteq U$ and $V_N \subseteq V$ be $N-$dimensional subspaces. We project the original problem onto $U_n$, i.e.,
 \begin{equation}
Find\quad u_n \in U_n\quad  s.t \quad \quad a(u_n,v_n)=\ell(v_n)  \quad \forall v_n \in V_n.
\label{eq:galerkin}
\end{equation}
Reducing the original problem to a $N-$dimensional subspace allows us to numerically solve for $u_n$ as a finite linear combination of basis vectors in $U_n$. For a finite number of basis functions, the Galerkin methods lead to a system of equations with a finite number of unknowns. Depending on the choice of test functions, these methods can be classified into different methods such as Bubnov-Galerkin, Petrov-Galerkin, collocation methods, or the least squares method. 
\par
Although the Galerkin methods can be used for solving a wide variety of problems, we restrict our attention to using ideas from Galerkin methods to develop time integrators for dynamical systems. Consider an autonomous system of first-order ODEs given by
\begin{equation}
\frac{\text{d}\y}{\text{d}t} = \f(\y), \quad \quad \y(t_0)=\y_0,
\end{equation}
 on a finite time interval $t \in [t_0,t_f]$ where $\y(t_0)=\y_0$ is the initial condition for the given initial value problem. We rewrite the governing first-order equations in the residual form by defining 
\begin{equation}
\R\left(\y,\frac{\text{d}\y}{\text{d}t}\right)= \frac{\text{d}\y}{\text{d}t} - \f(\y).
\end{equation}
We discretize the total time interval into $N$ subintervals of fixed time step size $\dt=\frac{t_f-t_0}{N}$. On each subinterval $[t_k, t_{k+1}]$, we consider piecewise smooth polynomial approximations of $\y$
\begin{equation}
\y(t) \approx \y_{\rm{d},k}(t) = \sum_{i=1}^{n+1} b^{(k)}_i \Phi_{i,k}(t),
\end{equation}
where $\y_{\rm{d},k}(t)$ is the approximation over the $k$th subinterval and  $\Phi_{i,k}(t)$ are basis functions that are $n$th degree polynomials on $k$th subinterval.  Since the coefficients $b_i^{(k)}$ may change from one subinterval to the next, the global approximation $\y_{\rm d}(t)$ need not be as smooth as the trial basis functions $\Phi_{i,k}(t)$. Using this piecewise polynomial approximation, we can derive a $\mathcal{C}^0-$continuous numerical integrator by requiring that $\y_{\rm{d},k}(t)$ is a local Galerkin approximation on the $k$th subinterval  and $\y_{\rm d}$ is continuous on $[t_0,t_f]$. For example, if we use $w_i=\Phi_{i,k}$ as test functions then that leads to solving the following $n+1$ equations for $k=0,\cdots,N-1$ 
\begin{equation}
\y_{\rm{d},k}(t_{k})=\begin{cases}\y_{\rm{d},k-1}(t_{k}), \quad & k\geq 1 \\ \y_0,\quad & k=0 \end{cases} \end{equation}
\begin{equation}
   \bigg \langle \textbf{R}\left(\y_{\rm d},\frac{\text{d}\y_{\rm d}}{\text{d}t}\right), \Phi_{i,k} \bigg \rangle_k = 0, \quad i=2, \cdots,n+1, 
\end{equation}
where $\langle v,w\rangle_k=\int_{t_k}^{t_{k+1}}v(t)w(t)\text{d}t$ is the inner product on the $k$th subinterval.}}
\section{One-step Time-integration Methods}
\label{s:3}
\subsection{Hermite Polynomials}
\label{s:Herm}
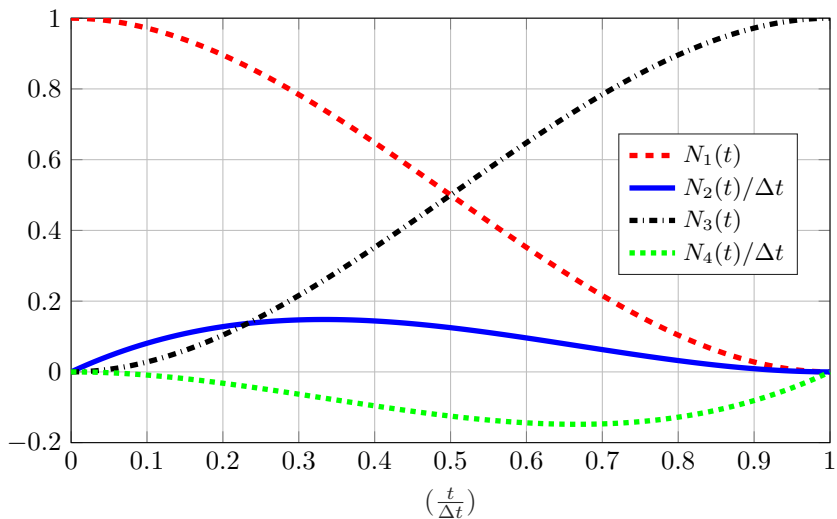
\begin{figure}[h]
\centering
    \setlength\fheight{10.5 cm}
        \setlength\fwidth{\textwidth}
%
%
\begin{tikzpicture}

\begin{axis}[%
width=0.951\fheight,
height=0.536\fheight,
at={(0\fheight,0\fheight)},
scale only axis,
xmin=0,
xmax=1,
xlabel style={font=\color{white!15!black}},
xlabel={ $(\frac{t}{\Delta t})$},
ymin=-0.2,
ymax=1,
axis background/.style={fill=white},
xmajorgrids,
ymajorgrids,
legend style={at={(0.721,0.395)}, anchor=south west, legend cell align=left, align=left, draw=white!15!black,font=\small }
]
\addplot [color=red, dashed, line width=2.0pt]
  table[row sep=crcr]{%
0	1\\
0.01	0.999702\\
0.02	0.998816\\
0.03	0.997354\\
0.04	0.995328\\
0.05	0.99275\\
0.0600000000000001	0.989632\\
0.0700000000000001	0.985986\\
0.0800000000000001	0.981824\\
0.0900000000000001	0.977158\\
0.1	0.972\\
0.11	0.966362\\
0.12	0.960256\\
0.13	0.953694\\
0.14	0.946688\\
0.15	0.93925\\
0.16	0.931392\\
0.17	0.923126\\
0.18	0.914464\\
0.19	0.905418\\
0.2	0.896\\
0.21	0.886222\\
0.22	0.876096\\
0.23	0.865634\\
0.24	0.854848\\
0.25	0.84375\\
0.26	0.832352\\
0.27	0.820666\\
0.28	0.808704\\
0.29	0.796478\\
0.3	0.784\\
0.31	0.771282\\
0.32	0.758336\\
0.33	0.745174\\
0.34	0.731808\\
0.35	0.71825\\
0.36	0.704512\\
0.37	0.690606\\
0.38	0.676544\\
0.39	0.662338\\
0.4	0.648\\
0.42	0.618976\\
0.44	0.589568\\
0.46	0.559872\\
0.49	0.514998\\
0.53	0.455054\\
0.55	0.42525\\
0.57	0.395686\\
0.59	0.366458\\
0.6	0.352\\
0.61	0.337662\\
0.62	0.323456\\
0.63	0.309394\\
0.64	0.295488\\
0.65	0.28175\\
0.66	0.268192\\
0.67	0.254826\\
0.68	0.241664\\
0.69	0.228718\\
0.7	0.216\\
0.71	0.203522\\
0.72	0.191296\\
0.73	0.179334\\
0.74	0.167648\\
0.75	0.15625\\
0.76	0.145152\\
0.77	0.134366\\
0.78	0.123904\\
0.79	0.113778\\
0.8	0.104\\
0.81	0.0945819999999999\\
0.82	0.0855360000000001\\
0.83	0.0768740000000001\\
0.84	0.068608\\
0.85	0.0607500000000001\\
0.86	0.053312\\
0.87	0.046306\\
0.88	0.039744\\
0.89	0.0336380000000001\\
0.9	0.028\\
0.91	0.022842\\
0.92	0.018176\\
0.93	0.014014\\
0.94	0.0103679999999999\\
0.95	0.00724999999999998\\
0.96	0.00467200000000001\\
0.97	0.00264599999999993\\
0.98	0.00118400000000007\\
0.99	0.000297999999999909\\
1	0\\
};
\addlegendentry{$N_1(t)$}

\addplot [color=blue, line width=2.0pt]
  table[row sep=crcr]{%
0	0\\
0.01	0.00980099999999995\\
0.02	0.0192079999999999\\
0.03	0.028227\\
0.04	0.036864\\
0.05	0.0451250000000001\\
0.0600000000000001	0.053016\\
0.0700000000000001	0.060543\\
0.0800000000000001	0.067712\\
0.0900000000000001	0.0745290000000001\\
0.1	0.081\\
0.11	0.0871310000000001\\
0.12	0.0929279999999999\\
0.13	0.0983970000000001\\
0.14	0.103544\\
0.15	0.108375\\
0.16	0.112896\\
0.17	0.117113\\
0.18	0.121032\\
0.19	0.124659\\
0.2	0.128\\
0.21	0.131061\\
0.22	0.133848\\
0.23	0.136367\\
0.24	0.138624\\
0.25	0.140625\\
0.26	0.142376\\
0.27	0.143883\\
0.28	0.145152\\
0.29	0.146189\\
0.3	0.147\\
0.31	0.147591\\
0.32	0.147968\\
0.33	0.148137\\
0.34	0.148104\\
0.35	0.147875\\
0.36	0.147456\\
0.37	0.146853\\
0.38	0.146072\\
0.39	0.145119\\
0.4	0.144\\
0.41	0.142721\\
0.42	0.141288\\
0.43	0.139707\\
0.44	0.137984\\
0.45	0.136125\\
0.46	0.134136\\
0.47	0.132023\\
0.48	0.129792\\
0.49	0.127449\\
0.5	0.125\\
0.51	0.122451\\
0.52	0.119808\\
0.53	0.117077\\
0.54	0.114264\\
0.55	0.111375\\
0.57	0.105393\\
0.59	0.0991789999999999\\
0.61	0.092781\\
0.63	0.086247\\
0.66	0.0762959999999999\\
0.7	0.0629999999999999\\
0.72	0.0564480000000001\\
0.74	0.0500240000000001\\
0.76	0.043776\\
0.78	0.037752\\
0.79	0.0348390000000001\\
0.8	0.032\\
0.81	0.0292410000000001\\
0.82	0.0265679999999999\\
0.83	0.023987\\
0.84	0.021504\\
0.85	0.0191250000000001\\
0.86	0.016856\\
0.87	0.0147029999999999\\
0.88	0.012672\\
0.89	0.010769\\
0.9	0.0089999999999999\\
0.91	0.00737100000000002\\
0.92	0.00588799999999989\\
0.93	0.00455699999999992\\
0.94	0.00338400000000005\\
0.95	0.00237500000000002\\
0.96	0.00153599999999998\\
0.97	0.000872999999999902\\
0.98	0.000391999999999948\\
0.99	9.90000000000713e-05\\
1	0\\
};
\addlegendentry{$N_2(t)/\Delta t$}

\addplot [color=black, dashdotted, line width=2.0pt]
  table[row sep=crcr]{%
0	0\\
0.01	0.000297999999999909\\
0.02	0.00118400000000007\\
0.03	0.00264599999999993\\
0.04	0.00467200000000001\\
0.05	0.00724999999999998\\
0.0600000000000001	0.0103679999999999\\
0.0700000000000001	0.014014\\
0.0800000000000001	0.018176\\
0.0900000000000001	0.022842\\
0.1	0.028\\
0.11	0.0336380000000001\\
0.12	0.039744\\
0.13	0.046306\\
0.14	0.053312\\
0.15	0.0607500000000001\\
0.16	0.068608\\
0.17	0.0768740000000001\\
0.18	0.0855360000000001\\
0.19	0.0945819999999999\\
0.2	0.104\\
0.21	0.113778\\
0.22	0.123904\\
0.23	0.134366\\
0.24	0.145152\\
0.25	0.15625\\
0.26	0.167648\\
0.27	0.179334\\
0.28	0.191296\\
0.29	0.203522\\
0.3	0.216\\
0.31	0.228718\\
0.32	0.241664\\
0.33	0.254826\\
0.34	0.268192\\
0.35	0.28175\\
0.36	0.295488\\
0.37	0.309394\\
0.38	0.323456\\
0.39	0.337662\\
0.4	0.352\\
0.42	0.381024\\
0.44	0.410432\\
0.46	0.440128\\
0.49	0.485002\\
0.53	0.544946\\
0.55	0.57475\\
0.57	0.604314\\
0.59	0.633542\\
0.6	0.648\\
0.61	0.662338\\
0.62	0.676544\\
0.63	0.690606\\
0.64	0.704512\\
0.65	0.71825\\
0.66	0.731808\\
0.67	0.745174\\
0.68	0.758336\\
0.69	0.771282\\
0.7	0.784\\
0.71	0.796478\\
0.72	0.808704\\
0.73	0.820666\\
0.74	0.832352\\
0.75	0.84375\\
0.76	0.854848\\
0.77	0.865634\\
0.78	0.876096\\
0.79	0.886222\\
0.8	0.896\\
0.81	0.905418\\
0.82	0.914464\\
0.83	0.923126\\
0.84	0.931392\\
0.85	0.93925\\
0.86	0.946688\\
0.87	0.953694\\
0.88	0.960256\\
0.89	0.966362\\
0.9	0.972\\
0.91	0.977158\\
0.92	0.981824\\
0.93	0.985986\\
0.94	0.989632\\
0.95	0.99275\\
0.96	0.995328\\
0.97	0.997354\\
0.98	0.998816\\
0.99	0.999702\\
1	1\\
};
\addlegendentry{$N_3(t)$}

\addplot [color=green, dotted, line width=2.0pt]
  table[row sep=crcr]{%
0	0\\
0.01	-9.90000000000713e-05\\
0.02	-0.000391999999999948\\
0.03	-0.000872999999999902\\
0.04	-0.00153599999999998\\
0.05	-0.00237500000000002\\
0.0600000000000001	-0.00338400000000005\\
0.0700000000000001	-0.00455699999999992\\
0.0800000000000001	-0.00588799999999989\\
0.0900000000000001	-0.00737100000000002\\
0.1	-0.0089999999999999\\
0.11	-0.010769\\
0.12	-0.012672\\
0.13	-0.0147029999999999\\
0.14	-0.016856\\
0.15	-0.0191250000000001\\
0.16	-0.021504\\
0.17	-0.023987\\
0.18	-0.0265679999999999\\
0.19	-0.0292410000000001\\
0.2	-0.032\\
0.21	-0.0348390000000001\\
0.22	-0.037752\\
0.24	-0.043776\\
0.26	-0.0500240000000001\\
0.28	-0.0564480000000001\\
0.3	-0.0629999999999999\\
0.34	-0.0762959999999999\\
0.37	-0.086247\\
0.39	-0.092781\\
0.41	-0.0991789999999999\\
0.43	-0.105393\\
0.45	-0.111375\\
0.46	-0.114264\\
0.47	-0.117077\\
0.48	-0.119808\\
0.49	-0.122451\\
0.5	-0.125\\
0.51	-0.127449\\
0.52	-0.129792\\
0.53	-0.132023\\
0.54	-0.134136\\
0.55	-0.136125\\
0.56	-0.137984\\
0.57	-0.139707\\
0.58	-0.141288\\
0.59	-0.142721\\
0.6	-0.144\\
0.61	-0.145119\\
0.62	-0.146072\\
0.63	-0.146853\\
0.64	-0.147456\\
0.65	-0.147875\\
0.66	-0.148104\\
0.67	-0.148137\\
0.68	-0.147968\\
0.69	-0.147591\\
0.7	-0.147\\
0.71	-0.146189\\
0.72	-0.145152\\
0.73	-0.143883\\
0.74	-0.142376\\
0.75	-0.140625\\
0.76	-0.138624\\
0.77	-0.136367\\
0.78	-0.133848\\
0.79	-0.131061\\
0.8	-0.128\\
0.81	-0.124659\\
0.82	-0.121032\\
0.83	-0.117113\\
0.84	-0.112896\\
0.85	-0.108375\\
0.86	-0.103544\\
0.87	-0.0983970000000001\\
0.88	-0.0929279999999999\\
0.89	-0.0871310000000001\\
0.9	-0.081\\
0.91	-0.0745290000000001\\
0.92	-0.067712\\
0.93	-0.060543\\
0.94	-0.053016\\
0.95	-0.0451250000000001\\
0.96	-0.036864\\
0.97	-0.028227\\
0.98	-0.0192079999999999\\
0.99	-0.00980099999999995\\
1	0\\
};
\addlegendentry{$N_4(t)/\Delta t$}

\end{axis}

\begin{axis}[%
width=1.227\fheight,
height=0.658\fheight,
at={(-0.16\fheight,-0.072\fheight)},
scale only axis,
xmin=0,
xmax=1,
ymin=0,
ymax=1,
axis line style={draw=none},
ticks=none,
axis x line*=bottom,
axis y line*=left
]
\end{axis}
\end{tikzpicture}%
 \caption{{\Ra{Cubic Hermite polynomials plotted over one time step.}}}
 \label{fig:Hermite}
\end{figure}

We discretize the continuous trajectory over one time step using cubic Hermite polynomials 
\begin{equation}
\q(t)\approx \q_{\rm d}(t) = \q_0N_1(t) + \textbf{v}_0N_2(t) +  \q_{1}N_3(t) + \textbf{v}_{1}N_4(t),
\label{eq:hermite}
\end{equation}
where $\q_{\rm d}(t)$ is the discrete approximation and the Hermite polynomials are given by
\begin{multline*}
N_1(t)=2 \left( \frac{t}{\dt} \right)^3 - 3\left( \frac{t}{\dt} \right)^2  + 1, \hfill  N_3(t)=-2 \left( \frac{t}{\dt} \right)^3 + 3\left( \frac{t}{\dt} \right)^2, \\  N_2(t) = \dt \left[  \left( \frac{t}{\dt} \right)^3 - 2\left( \frac{t}{\dt}   \right)^2 + \left( \frac{t}{\dt}   \right)  \right], \hfill  N_4(t) = \dt \left[  \left( \frac{t}{\dt} \right)^3 - \left( \frac{t}{\dt}   \right)^2  \right],
\end{multline*}
where $\dt$ is the fixed time step. As shown in Figure \ref{fig:Hermite}, at initial time $t=0$, we have $N_1(0)=1$ along with $N_2(0)=N_3(0)=N_4(0)=0$ which leads to $\q_{\rm d}(0)=\q_0$. Similarly, we also have $\dot{\q}_{\rm d}(0)=\textbf{v}_0$, $\q_{\rm d}(\dt)=\q_1$ and $\dot{\q}_{\rm d}(\dt)=\textbf{v}_1$. Using piecewise Hermite polynomials for discretization in the one-step approach leads to numerical solutions that are $\mathcal{C}^1$-- continuous.
\subsection{One-step Variational Methods} 
\label{s:ov}
Unlike the discrete mechanics framework developed by Marsden and West \cite{marsden2001discrete}, we only consider the variational principle over one time step $\dt$ and derive one-step numerical integrators from the discretized variational principle.
\par
{\Rb{We use the discrete approximation $\q_{\rm d}(t)$ from \eqref{eq:hermite} to obtain a discrete action $S_{\rm d}$ which approximates the action integral over one time step in the following sense 
\begin{equation}
    S_{\rm d}(\q_0,\textbf{v}_0,\q_{1},\textbf{v}_{1}) = \int_{0}^{\dt} L(\q_{\rm d},\dot{\q}_{\rm d})  \  \text{d}t \approx \int_{0}^{\dt} L(\q,\dot{\q})  \  \text{d}t.
\end{equation}
Thus, we have used a cubic Hermite polynomial to obtain the discrete action $S_{\rm d}$ in terms of the configuration and velocity at the endpoints.}} We introduce discrete forces corresponding to displacement and velocity variations to approximate the virtual work
\begin{multline}
\int ^{\dt}_{0} \textbf{f}_L (\q(t), \dot{\q}(t)) \cdot\delta \q \ \text{d}t \approx \int  ^{\dt}_{0} \textbf{f}_L (\q_{\rm d}(t), \dot{\q}_{\rm d}(t)) \cdot\delta \q_{\rm d} \ \text{d}t \\
=\textbf{f}_{\rm d}^+(\q_0,\textbf{v}_0,\q_{1},\textbf{v}_{1})\cdot\delta \q_{1} +  \textbf{f}_{\rm d}^-(\q_0,\textbf{v}_0,\q_{1},\textbf{v}_{1})\cdot\delta \q_{0} \ \ \ \ \ \ \ \ \  \ \ \ \ \   \\
+ \textbf{g}_{\rm d}^+(\q_0,\textbf{v}_0,\q_{1},\textbf{v}_{1})\cdot\delta \textbf{v}_{1} +  \textbf{g}_{\rm d}^-(\q_0,\textbf{v}_0,\q_{1},\textbf{v}_{1})\cdot\delta \textbf{v}_{0}. 
\end{multline}
Using $\delta \q_{\rm d}(t) = \delta \q_0N_1(t) + \delta \textbf{v}_0N_2(t) +  \delta \q_{1}N_3(t) + \delta \textbf{v}_{1}N_4(t)$, we can obtain the discrete forces $\textbf{f}_{\rm d}^{\pm}$ corresponding to displacement variations
\begin{equation}
\textbf{f}_{\rm d}^+(\q_0,\textbf{v}_0,\q_{1},\textbf{v}_{1}) = \int  ^{\dt}_{0} \textbf{f}_L (\q_{\rm d}(t), \dot{\q}_{\rm d}(t)) \ N_3(t) \ \text{d}t,
\end{equation}
\begin{equation}
\textbf{f}_{\rm d}^-(\q_0,\textbf{v}_0,\q_{1},\textbf{v}_{1}) = \int  ^{\dt}_{0} \textbf{f}_L (\q_{\rm d}(t), \dot{\q}_{\rm d}(t)) \ N_1(t) \ \text{d}t,
\end{equation}
 and discrete forces $\textbf{g}_{\rm d}^{\pm}$ corresponding to velocity variations
\begin{equation}
\textbf{g}_{\rm d}^+(\q_0,\textbf{v}_0,\q_{1},\textbf{v}_{1}) = \int  ^{\dt}_{0} \textbf{f}_L (\q_{\rm d}(t), \dot{\q}_{\rm d}(t)) \ N_4(t) \ \text{d}t,
\end{equation}
\begin{equation}
\textbf{g}_{\rm d}^-(\q_0,\textbf{v}_0,\q_{1},\textbf{v}_{1}) = \int  ^{\dt}_{0} \textbf{f}_L (\q_{\rm d}(t), \dot{\q}_{\rm d}(t)) \ N_2(t) \ \text{d}t.
\end{equation}
{\Rb{The discrete Lagrange-d'Alembert principle  using the one-step variational approach seeks $ \q_{\rm d}(t)$  \eqref{eq:hermite} that satisfy 
\begin{multline*}
\delta S_{\rm d}(\q_0,\textbf{v}_0, \q_{1},\textbf{v}_{1}) + [\textbf{f}_{\rm d}^+(\q_0,\textbf{v}_0, \q_{1},\textbf{v}_{1})\cdot\delta \q_{1} +  \textbf{f}_{\rm d}^-(\q_0,\textbf{v}_0, \q_{1},\textbf{v}_{1})\cdot\delta \q_{0}] \\ + [\textbf{g}_{\rm d}^+(\q_0,\textbf{v}_0, \q_{1},\textbf{v}_{1})\cdot\delta \textbf{v}_{1} +  \textbf{g}_{\rm d}^-(\q_0,\textbf{v}_0, \q_{1},\textbf{v}_{1})\cdot\delta \textbf{v}_{0}]  = 0,
\end{multline*}
which gives
\begin{multline*}
 \left(  \frac{\partial S_{\rm d}}{\partial \q_0} + \textbf{f}_{\rm d}^-(\q_0,\textbf{v}_0, \q_{1},\textbf{v}_{1}) \right) \cdot \delta \q_0 + \left(  \frac{\partial S_{\rm d}}{\partial \textbf{v}_0} + \textbf{g}_{\rm d}^-(\q_0,\textbf{v}_0, \q_{1},\textbf{v}_{1}) \right) \cdot \delta \textbf{v}_0 \\ +  \left(  \frac{\partial S_{\rm d}}{\partial \q_1} + \textbf{f}_{\rm d}^+(\q_0,\textbf{v}_0, \q_{1},\textbf{v}_{1}) \right) \cdot \delta \q_1 + \left(  \frac{\partial S_{\rm d}}{\partial \textbf{v}_1} + \textbf{g}_{\rm d}^+(\q_0,\textbf{v}_0, \q_{1},\textbf{v}_{1}) \right) \cdot \delta \textbf{v}_1  =0,
\end{multline*}
where $S_{\rm d}:=S_{\rm d}(\q_0,\textbf{v}_0, \q_{1},\textbf{v}_{1}).$ Setting variations at endpoints to zero, i.e. $\delta \q_0= \delta \q_1=\textbf{0}$, gives 
\begin{equation}
 \frac{\partial S_{\rm d}}{\partial \textbf{v}_0}(\q_0,\textbf{v}_0, \q_{1},\textbf{v}_{1}) + \textbf{g}_{\rm d}^-(\q_0,\textbf{v}_0, \q_{1},\textbf{v}_{1})=\textbf{0},
 \end{equation}
 \begin{equation}
  \frac{\partial S_{\rm d}}{\partial \textbf{v}_1}(\q_0,\textbf{v}_0, \q_{1},\textbf{v}_{1}) + \textbf{g}_{\rm d}^+(\q_0,\textbf{v}_0, \q_{1},\textbf{v}_{1}) =\textbf{0}.
 \end{equation}
 Given $(\q_0,\textbf{v}_0)$, these coupled nonlinear equations can be solved to obtain $(\q_1,\textbf{v}_1)$.}} Thus, for Lagrangian systems with external forcing the one-step variational approach can be used to derive numerical integrators that are continuous in both configuration and velocities.
 {\Rb{
 \begin{remark}
It is important to note that we have used the time finite elements approach \cite{argyris1969finite} where considering the discrete Lagrange-d'Alembert principle over a single time step plays a crucial role in ensuring that the numerical integrators are $\mathcal{C}^1-$continuous and stable. Riff and Baruch \cite{riff1984stability} have shown that numerical integrators derived by considering the discrete action sum over the entire time interval are unconditionally unstable.
 \end{remark}
  \begin{remark}
In \cite{leok2012prolongation}, Hermite-based prolongation-collocation variational integrators (PCVIs) are constructed from discrete Lagrangian $L_{\rm d}(\q_0,\q_1)$ by means of expressing every parameter in $\q_{\rm d}(t)$ as function of $(\q_0,\q_1)$ and with the help of some extra equations based on the prolongation-collocation approach. PCVIs are fundamentally different from the one-step variational methods developed here in two ways. First, one-step variational methods proposed in this work are $\mathcal{C}^1-$continuous whereas the PCVIs in \cite{leok2012prolongation} are $\mathcal{C}^0-$continuous. Second, PCVIs are constructed by summing up the discrete action over the entire time interval $T$ whereas the one-step variational methods are constructed by considering the discrete action $S_{\rm d}(\q_0,\textbf{v}_0, \q_{1},\textbf{v}_{1})$ over a single time step $\Delta t$.   \end{remark}
}}
\subsection{One-step Galerkin Methods}
In this subsection, we consider the Petrov-Galerkin method or the weighted residual method to obtain numerical integrators for mechanical systems with external forcing. For a given system of equations we first write it in the residual form $\textbf{R}\left(t,\q_{\rm d},\frac{\text{d}\q_{\rm d}}{\text{d}t},...,\frac{\text{d}^n\q_{\rm d}}{\text{d}t^n}\right)$ and this continuous system of ODEs is transformed into following discrete equations
\begin{equation}
   \bigg \langle \textbf{R}\left(t,\q_{\rm d},\frac{\text{d}\q_{\rm d}}{\text{d}t},...,\frac{\text{d}^n\q_{\rm d}}{\text{d}t^n}\right), w_i \bigg \rangle = 0, \quad i=1,2,...,N,
\end{equation}
where $\q_{\rm d}(t)$ is the assumed solution form and  $w_i$ are test functions. For our study, we focus on time-integration of mechanical problems with $\textbf{R}=\textbf{M}_{\rm eq}(\q,\dot{\q},\ddot{\q})$ where $\textbf{M}_{\rm eq}$ are the equations of motion for the mechanical system. Similar to the variational approach discussed in the previous section, we use cubic Hermite polynomials as solution functions for approximating the continuous solution over one time step. Our goal is to use the Petrov-Galerkin approach to derive one-step methods so we consider the following shifted Legendre polynomials as test functions
\begin{equation*}
    P_0(t) =1, \quad P_1(t)= \frac{1}{\dt} (2t-\dt),
\end{equation*}
Given $(\q_0,\textbf{v}_0)$, the one-step Galerkin method yields 
\begin{equation*}
\int_0^{\dt} \textbf{M}_{\rm eq}(\q_{\rm d}(t),\dot{\q}_{\rm d}(t),\ddot{\q}_{\rm d}(t))\left( 1 \right) \text{d}t=\textbf{0},  
\end{equation*}
\begin{equation*}
\int_0^{\dt} \textbf{M}_{\rm eq}(\q_{\rm d}(t),\dot{\q}_{\rm d}(t),\ddot{\q}_{\rm d}(t))\left(\frac{1}{\dt}(2t-\dt)\right)\text{d}t=\textbf{0}. 
\end{equation*}
Given  $(\q_0,\textbf{v}_0)$, these two coupled nonlinear equations can be solved to obtain the configuration $\q_1$ and velocity $\textbf{v}_1$ at the next time step. Just like the variational approach, this system of nonlinear equations can be used as a one-step numerical integrator.
\subsection{One-step Methods based on Higher-order Hermite Polynomials}
In this subsection, we demonstrate how to derive one-step methods with higher-order Hermite polynomials. The one-step methods presented so far have been based on cubic Hermite polynomials which lead to numerical integrators that are continuous in both configuration and velocities. For discretization using higher-order Hermite polynomials, we consider the following discrete trajectory over one time step
\begin{equation}
    \q_{\rm d}(t)=\sum_{j=0}^{n-1} \left( \q^{(j)}_0 H_{n,j}(t) + \q^{(j)}_1 H_{n,j}(\dt-t) \right),
\end{equation}
where we have written the discrete trajectory in terms of Hermite basis functions and values of $\q(t)$ and its derivatives $\q^{(j)}(t)$ at endpoints of each interval. The Hermite basis functions are
\begin{equation}
    H_{n,j}(t)= \frac{t^j}{j!} \left(1-\frac{t}{\dt}\right)^n \sum_{s=0}^{n-j-1} \begin{pmatrix}
n+s-1 \\
s
\end{pmatrix}\left( \frac{t}{\dt} \right)^s.
\end{equation}
Thus, the discrete trajectory $\q_{\rm d}$ is represented by a $2n-1$ degree polynomial which satisfies 
\begin{equation}
  \q^{(j)}(0)=\q_{\rm d}^{(j)}(0), \quad \quad \quad  \q^{(j)}(\dt)=\q_{\rm d}^{(j)}(\dt),   \quad \quad j=0,...,n-1.
\end{equation}
It is clear from the above expression that for  $n=1$, the discrete trajectory simply reduces to a linear interpolant between endpoints $\q_0$ and $\q_1$. For $n=2$, the discrete trajectory is the cubic interpolant discussed in Section \ref{s:Herm}. For discretization using higher-order Hermite polynomials with $n \geq 3$, the discrete trajectory $\q_{\rm d}$ over the fixed time step is represented by $2n-1$ degree polynomials with $2n$ unknown coefficients. For an initial condition of the form $(\q_0,\textbf{v}_0)$, the numerical integration problem reduces to solving for the remaining $2n-2$ coefficients. The first $n-2$ coefficients are $\q^{j}(0)$ for $j=2,\cdots, n-1$ and the other $n$ coefficients are $\q^{j}(\dt)$ for $j=0,\cdots, n-1$. 
\\ \\
\textbf{Variational Approach:} For conservative Lagrangian systems, the discrete Hamilton's principle after setting the configuration variations at the endpoints to zero (i.e. $\delta \q_0=\delta \q_1=\textbf{0}$) leads to the following discrete equations
\begin{equation}
    \frac{\partial S_{\rm d}}{\partial \q^{(j)}(0)}= \frac{\partial S_{\rm d}}{\partial \q^{(j)}(\dt)} = \textbf{0}, \quad \quad j=1,\cdots, n-1.
\end{equation}
Thus, solving these $2n-2$ coupled nonlinear equations gives the $2n-2$ coefficients and this one-step method can be seen as a numerical integrator from $(\q_0,\textbf{v}_0)$ to $(\q_1,\textbf{v}_1)$ with a $2n-1$ degree Hermite piecewise polynomial interpolating the configuration over every fixed time step. This approach can be extended to Lagrangian systems with forcing by discretizing the Lagrange-d'Alembert principle. 
\\ \\
\textbf{Galerkin Approach:} One-step Galerkin methods with discretization using higher-order Hermite polynomials involve the use of shifted Legendre polynomials up to order $2n-2$ as test functions. Given $(\q_0,\textbf{v}_0)$, the governing discrete equations are given by
\begin{equation}
    \int_0^{\dt} \textbf{M}_{\rm eq}(\q_{\rm d}(t),\dot{\q}_{\rm d}(t),\ddot{\q}_{\rm d}(t)) P_j(t)  \ \text{d}t = \textbf{0}, \quad \quad j=0,\cdots, 2n-3,
\end{equation}
where $P_j(t)$ are shifted Legendre polynomials of degree $j$.
\section{Numerical Properties}
\label{s:4}
In order to understand the numerical properties of the proposed one-step methods, we investigate the linear stability and symplectic nature of these algorithms. We consider the simple harmonic oscillator with a single degree of freedom for both studies. 
\subsection{Linear Stability}
\begin{figure}
\captionsetup[subfigure]{oneside,margin={1.8cm,0 cm}}
    \begin{subfigure}{.37\textwidth}
       \setlength\fheight{5.5 cm}
        \setlength\fwidth{\textwidth}
\input{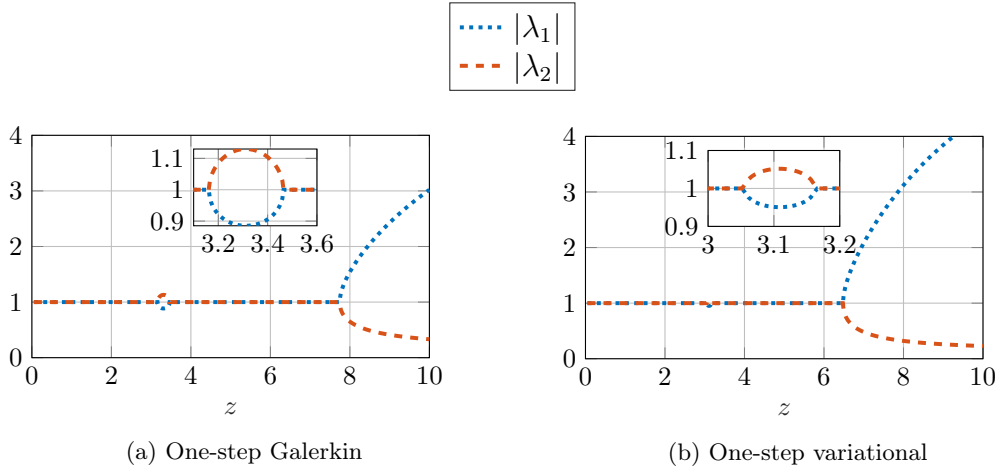}
 \caption{One-step Galerkin}
\label{fig:stab_gal}
    \end{subfigure}
    \hspace{1.5cm}
    \begin{subfigure}{.37\textwidth}
           \setlength\fheight{5.5 cm}
           \setlength\fwidth{\textwidth}
\raisebox{-5.6cm}{
%
%
\definecolor{mycolor1}{rgb}{0.00000,0.44700,0.74100}%
\definecolor{mycolor2}{rgb}{0.85000,0.32500,0.09800}%
\begin{tikzpicture}

\begin{axis}[%
width=0.951\fheight,
height=0.536\fheight,
at={(0\fheight,0\fheight)},
scale only axis,
xmin=0,
xmax=10,
xlabel style={font=\color{white!15!black}},
xlabel={$z$},
ymin=0,
ymax=4,
axis background/.style={fill=white},
xmajorgrids,
ymajorgrids,
legend style={at={(0.246,0.757)}, anchor=south west, legend cell align=left, align=left, draw=white!15!black}
]
\addplot [color=mycolor1, dotted, line width=1.5pt]
  table[row sep=crcr]{%
0.0500000000000007	1\\
3.05142571285643	1\\
3.06138069034517	0.976356598218386\\
3.06635817908954	0.96930987624429\\
3.07133566783392	0.96424252528805\\
3.07631315657829	0.960367663229851\\
3.08129064532266	0.957346547017144\\
3.0912456228114	0.953215092115057\\
3.10120060030015	0.951098570314837\\
3.11115557778889	0.950712754780771\\
3.12111055527764	0.952006099093513\\
3.13106553276638	0.955127870256298\\
3.13604302151075	0.957497137043214\\
3.14102051025513	0.960533765167826\\
3.1459979989995	0.964416660427339\\
3.15097548774387	0.969481387586635\\
3.15595297648824	0.976508676533319\\
3.16590795397699	1\\
6.47593796898449	1\\
6.48091545772886	1.0143225953936\\
6.48589294647324	1.08022926575059\\
6.49087043521761	1.11421494478539\\
6.49584792396198	1.14111017808083\\
6.50580290145072	1.18510386251368\\
6.51575787893947	1.22204276994373\\
6.52571285642821	1.25482447190489\\
6.53566783391696	1.28478061162379\\
6.55060030015007	1.32596252955347\\
6.56553276638319	1.36386308299638\\
6.58046523261631	1.39931883849993\\
6.6003751875938	1.44368088384959\\
6.62028514257129	1.48544100977526\\
6.64517258629315	1.53477875088105\\
6.67006003001501	1.58157392145684\\
6.69494747373687	1.62630081021191\\
6.7248124062031	1.67771304480181\\
6.75467733866934	1.7270548727419\\
6.78951975987994	1.78241330723748\\
6.82436218109054	1.83573956891174\\
6.85920460230115	1.88730790246588\\
6.89902451225613	1.94435975882139\\
6.9388444222111	1.99963651634707\\
6.98364182091045	2.05993976978912\\
7.0284392196098	2.11845724250971\\
7.07323661830915	2.17536242788698\\
7.12301150575288	2.23687010014624\\
7.1727863931966	2.29671926701288\\
7.22256128064032	2.35504043658365\\
7.27731365682842	2.4175586560464\\
7.33206603301651	2.47848208788325\\
7.3868184092046	2.5379142714541\\
7.44654827413707	2.60115259627558\\
7.50627813906953	2.66282103114265\\
7.566008004002	2.72300573202774\\
7.62573786893447	2.78178220867772\\
7.69044522261131	2.8439449318175\\
7.75515257628814	2.90460888023038\\
7.81985992996498	2.9638415933273\\
7.88456728364182	3.02170392535374\\
7.95425212606303	3.0825479210134\\
8.02393696848424	3.14192951628248\\
8.09362181090545	3.19990624370061\\
8.16330665332666	3.25653105980009\\
8.23299149574787	3.31185298925532\\
8.30765382691346	3.36973244427804\\
8.38231615807904	3.4262208944097\\
8.45697848924462	3.48136724543822\\
8.5316408204102	3.53521738957756\\
8.60630315157579	3.58781454734713\\
8.68594297148574	3.64258308111371\\
8.7655827913957	3.69601930537149\\
8.84522261130565	3.74816760157624\\
8.92486243121561	3.79907009312357\\
9.00450225112556	3.84876684828918\\
9.08911955977989	3.90029120465912\\
9.17373686843422	3.95054115634234\\
9.25835417708855	3.99955834460584\\
9.34297148574287	4.04738253644675\\
9.4275887943972	4.09405175949785\\
9.5171835917959	4.14224779060554\\
9.6067783891946	4.1892311468719\\
9.6963731865933	4.23504141900257\\
9.78596798399199	4.27971653038672\\
9.88054027013507	4.32568223293765\\
9.97511255627814	4.37046461204006\\
10	4.3820577342606\\
};

\addplot [color=mycolor2, dashed, line width=1.5pt]
  table[row sep=crcr]{%
0.0500000000000007	1\\
3.05142571285643	1\\
3.0564032016008	1.01139856511699\\
3.06138069034517	1.02421594919802\\
3.06635817908954	1.03166182921258\\
3.07133566783392	1.03708348654429\\
3.07631315657829	1.04126787925872\\
3.08129064532266	1.04455382757242\\
3.08626813406703	1.04712232458016\\
3.09622311155578	1.05049803667618\\
3.10617808904452	1.05185980023253\\
3.11613306653327	1.0513640389084\\
3.12608804402201	1.04896606575241\\
3.13106553276638	1.04698023284742\\
3.13604302151075	1.04438954573591\\
3.14102051025513	1.04108781623651\\
3.1459979989995	1.03689623067782\\
3.15097548774387	1.03147931750329\\
3.15595297648824	1.02405644110616\\
3.16093046523262	1.01129823912779\\
3.16590795397699	1\\
6.47593796898449	1\\
6.48091545772886	0.985879644741585\\
6.48589294647324	0.925729409214954\\
6.49087043521761	0.897492898188162\\
6.49584792396198	0.876339567562042\\
6.50082541270635	0.858887367693454\\
6.50580290145072	0.843807898726226\\
6.51575787893947	0.818301965033553\\
6.52571285642821	0.796924209233779\\
6.53566783391696	0.778343003430081\\
6.5456228114057	0.761815130062759\\
6.55557778889445	0.746875783787052\\
6.57051025512756	0.72678507901937\\
6.58544272136068	0.70887220962981\\
6.6003751875938	0.692673852779352\\
6.61530765382691	0.677868549281373\\
6.6352176088044	0.659899807100937\\
6.65512756378189	0.643595259365632\\
6.67503751875938	0.62866497281718\\
6.69494747373687	0.61489239488831\\
6.71983491745873	0.599054246510679\\
6.74472236118059	0.584516178126202\\
6.76960980490245	0.571086022346243\\
6.79947473736868	0.556221991430217\\
6.82933966983492	0.542532192948675\\
6.86418209104552	0.5278314277629\\
6.89902451225613	0.514308113744427\\
6.9388444222111	0.500090887431282\\
6.97866433216608	0.487013506519052\\
7.02346173086543	0.473476379271871\\
7.06825912956478	0.461017590684074\\
7.1180340170085	0.448270415400984\\
7.1727863931966	0.435403670950436\\
7.22753876938469	0.423586504043126\\
7.28726863431716	0.411734260627108\\
7.351975987994	0.399958183530984\\
7.42166083041521	0.388349597539095\\
7.49632316158079	0.376981767909548\\
7.57596298149075	0.3659119766661\\
7.66058029014507	0.355183632077249\\
7.75515257628814	0.344280431973576\\
7.85470235117559	0.333871436439759\\
7.96420710355178	0.323507426649114\\
8.08366683341671	0.31331188763858\\
8.20810405202601	0.303744397989975\\
8.34747373686843	0.294115092131085\\
8.4967983991996	0.284880333510426\\
8.66105552776388	0.275815413572364\\
8.83526763381691	0.26725690976453\\
9.02938969484742	0.258795427941008\\
9.23844422211105	0.250743852906554\\
9.46740870435218	0.24297648902226\\
9.72126063031516	0.235429090126434\\
10	0.228203291842007\\
};

\end{axis}

\begin{axis}[%
width=0.315\fheight,
height=0.183\fheight,
at={(0.2933\fheight,0.319\fheight)},
scale only axis,
xmin=3,
xmax=3.2,
xtick={3,3.1,3.2},
ymin=0.9,
ymax=1.1,
ytick={0.9,1,1.1},
axis background/.style={fill=white},
xmajorgrids,
ymajorgrids
]
\addplot [color=mycolor1, dotted, line width=1.5pt, forget plot]
  table[row sep=crcr]{%
2.99667333666833	1\\
3.05142571285643	1\\
3.0564032016008	0.988729897875951\\
3.06138069034517	0.976356598218385\\
3.06635817908954	0.96930987624429\\
3.07133566783392	0.964242525288049\\
3.07631315657829	0.960367663229852\\
3.08129064532266	0.957346547017144\\
3.08626813406703	0.954998261927947\\
3.09124562281141	0.953215092115058\\
3.09622311155578	0.951929432599459\\
3.10120060030015	0.951098570314837\\
3.10617808904452	0.95069704135374\\
3.11115557778889	0.950712754780771\\
3.11613306653327	0.951145334054101\\
3.12111055527764	0.952006099093514\\
3.12608804402201	0.953319685592223\\
3.13106553276638	0.955127870256298\\
3.13604302151076	0.957497137043214\\
3.14102051025513	0.960533765167826\\
3.1459979989995	0.964416660427339\\
3.15097548774387	0.969481387586635\\
3.15595297648824	0.97650867653332\\
3.16093046523262	0.988827984969558\\
3.16590795397699	1\\
3.20075037518759	1\\
};
\addplot [color=mycolor2, dashed, line width=1.5pt, forget plot]
  table[row sep=crcr]{%
2.99667333666833	1\\
3.05142571285643	1\\
3.0564032016008	1.01139856511699\\
3.06138069034517	1.02421594919803\\
3.06635817908954	1.03166182921258\\
3.07133566783392	1.03708348654429\\
3.07631315657829	1.04126787925872\\
3.08129064532266	1.04455382757242\\
3.08626813406703	1.04712232458016\\
3.09124562281141	1.04908116570116\\
3.09622311155578	1.05049803667618\\
3.10120060030015	1.05141573251338\\
3.10617808904452	1.05185980023253\\
3.11115557778889	1.05184241504217\\
3.11613306653327	1.0513640389084\\
3.12111055527764	1.05041343847711\\
3.12608804402201	1.04896606575241\\
3.13106553276638	1.04698023284742\\
3.13604302151076	1.04438954573591\\
3.14102051025513	1.04108781623651\\
3.1459979989995	1.03689623067782\\
3.15097548774387	1.03147931750329\\
3.15595297648824	1.02405644110616\\
3.16093046523262	1.01129823912779\\
3.16590795397699	1\\
3.20075037518759	1\\
};
\end{axis}

\begin{axis}[%
width=1.227\fheight,
height=0.658\fheight,
at={(-0.16\fheight,-0.072\fheight)},
scale only axis,
xmin=0,
xmax=1,
ymin=0,
ymax=1,
axis line style={draw=none},
ticks=none,
axis x line*=bottom,
axis y line*=left
]
\end{axis}
\end{tikzpicture}%
}
 \caption{One-step variational}
\label{fig:stab_var}
    \end{subfigure}
    \caption{{\Ra{Linear stability analysis of both one-step methods.}}}
\label{fig:stability}
\end{figure}

We closely follow Leimkuhler and Reich \cite{leimkuhler_reich_2005} to study the linear stability of the proposed one-step methods for the simple harmonic oscillator. We consider the following Lagrangian system
\begin{equation}
    L(q,\dot{q})=\frac{1}{2}\dot{q}^2 - \frac{1}{2}\omega^2q^2,
\end{equation}
where $\omega$ is the natural frequency. The governing equation is
\begin{equation}
    \ddot{q} + \omega^2q^2 = 0.
\end{equation}
The analytical solution for this Lagrangian system is given by
\begin{equation}
 \begin{bmatrix}
    v(t) \\ \omega q(t) \end{bmatrix} = 
   \underbrace{ \begin{bmatrix}
    \cos(\omega t) & -\sin(\omega t) \\ 
    \sin(\omega t) & \cos(\omega t) \end{bmatrix}}_{=A_{\omega}}  
    \begin{bmatrix}
    v(0) \\ \omega q(0) \end{bmatrix} .
\end{equation}
Since $A_{\omega}A_{\omega}^\top =\mathbb{I}$, $A_{\omega}$ is orthogonal and thus, the eigenvalues are $\lambda_{1,2}=e^{\pm i \omega t}$ with $|\lambda_{1,2}|=1$.
For this linear dynamical system, time-marching equations for both one-step methods can be written in the following form
\begin{equation}
    \begin{bmatrix}
    v_{k+1} \\ \omega q_{k+1}
    \end{bmatrix} = {\Ra{A_{z} }} \begin{bmatrix}
    v_{k} \\ \omega q_k
    \end{bmatrix},
\end{equation}
where {\Ra{$A_{z}$ is the amplification matrix with $z:=\dt\omega$}}.  A sufficient condition for the asymptotic stability of a numerical method is that the eigenvalues of the amplification matrix {\Ra{$A_{z}$}} are in the unit disk of the complex plane and are simple if they lie on the unit circle. We investigate the linear stability of both Hermite-based one-step methods and compare the results with the variational integrators based on the discrete mechanics. {\Ra{
\begin{enumerate}
    \item For the one-step variational method we have
\begin{equation*}
     A_{z} = \begin{bmatrix}
  \frac{7z^4 - 192z^2 + 420 }{2z^4 + 18z^2 + 420}&    \frac{15z(3z^2 - 28)}{2z^4 + 18z^2 + 420} \\   \frac{z(z^4 - 52z^2 + 420) }{2z^4 + 18z^2 + 420} &    \frac{3z^4 - 104z^2 + 240 }{2z^4 + 18z^2 + 420} \end{bmatrix},
\end{equation*}
with eigenvalues $\lambda_{1,2}=\frac{17z^4-192z^2 + 420 \pm z \sqrt{15(z^2-10)(3z^2-28)(z^2-42)}}{2z^4 + 18z^2 + 420}$. The stability region for the one-step variational method is shown in Figure \ref{fig:stab_var}. We observe a small region of instability for $\sqrt{\frac{28}{3}}< z< \sqrt{10}$. Thus, the one-step variational method is stable for $z<\sqrt{\frac{28}{3}}$.
\item For the one-step Galerkin method we have
\begin{equation*}
     A_{z} = \begin{bmatrix}
  \frac{3z^4 - 104z^2 + 240 }{3z^4 + 48z^2 + 720}&    \frac{24z(z^2 - 10)}{3z^4 + 48z^2 + 720} \\   \frac{z(z^4 - 72z^2 + 720) }{3z^4 + 48z^2 + 720} &    \frac{3z^4 - 104z^2 + 240 }{3z^4 + 48z^2 + 720} \end{bmatrix},
\end{equation*}
with eigenvalues $\lambda_{1,2}=\frac{3z^4-104z^2 + 240 \pm 2z \sqrt{2(z^2-10)(z^2-12)(z^2-60)}}{z^4 + 16z^2 + 240}$. The stability region for the one-step Galerkin method is shown in Figure \ref{fig:stab_gal}. Simillar to the one-step variational method, we observe a small region of instability for $\sqrt{10}< z< \sqrt{12}$. Thus, the one-step Galerkin method is stable for $z<\sqrt{10}$.
\item The amplification matrix $A_{z}$ for the variational integrator based on the midpoint rule is 
\begin{equation*}
     A_{z} = \begin{bmatrix}
  \frac{4-z^2 }{ z^2 + 4}&    \frac{4z}{z^2 + 4} \\   \frac{-4z }{z^2 + 4} &    \frac{4-z^2  }{z^2 + 4} \end{bmatrix},
\end{equation*}
The amplification matrix for this method is orthogonal and hence, the method is stable for all $z \in \mathbb{R}$.
\end{enumerate}
The eigenvalues for the proposed one-step methods are plotted for different $z$ values in Figure \ref{fig:stability} where both one-step methods have similar stability regions. The one-step Galerkin method is stable for $\dt<\frac{\sqrt{10}}{\omega}$ whereas the one-step variational method is stable for $\dt < \frac{\sqrt{\frac{28}{3}}}{\omega}$. }}
\subsection{Symplectic Nature}
As mentioned in Section \ref{s:ov}, the variational approach to one-step methods only considers the action integral over one time step whereas the variational integrators consider the action integral over a finite number of fixed time steps. We know that variational integrators derived from the latter approach yield numerical algorithms that automatically preserve the canonical symplectic form. For a Hamiltonian system with Hamiltonian $H(\textbf{p},\q)$, the  symplectic flow map $\phi_t(\textbf{p}_0,\q_0)=(\textbf{p}(t),\q(t))$ satisfies the following condition
\begin{equation}
    \left( \frac{\partial \phi_t}{\partial \textbf{y}_0} \right)^\top J \left( \frac{\partial \phi_t}{\partial \textbf{y}_0} \right) =J,
\end{equation}
where $\textbf{y}_0=(\textbf{p}_0,\q_0)$ and $J=\begin{bmatrix}
0 & I \\ -I & 0
\end{bmatrix}$ is the symplectic matrix. Similar to this condition, a given one-step method $\phi_{\dt}:(\textbf{p}_k,\q_k) \to (\textbf{p}_{k+1},\q_{k+1})$ is symplectic if it satisfies $ \left( \frac{\partial \phi_{\dt}}{\partial \textbf{y}_k} \right)^\top J \left( \frac{\partial \phi_{\dt}}{\partial \textbf{y}_k} \right) =J$ for $\textbf{y}_k=(\textbf{p}_k,\q_k)$. The key step in this process is to compute the following Jacobian matrix 
\begin{equation}
    \frac{\partial \phi_{\dt}}{\partial \textbf{y}_k} = \begin{bmatrix}
  \frac{\partial \textbf{p}_{k+1}(\textbf{p}_k,\q_k) }{\partial \textbf{p}_k}&  \frac{\partial \textbf{p}_{k+1}(\textbf{p}_k,\q_k) }{\partial \q_k} \\ \frac{\partial \q_{k+1}(\textbf{p}_k,\q_k) }{\partial \textbf{p}_k}&  \frac{\partial \q_{k+1}(\textbf{p}_k,\q_k) }{\partial \q_k}
\end{bmatrix}.
\end{equation}
Since both proposed methods are generally implicit, the computation for the Jacobian matrix involves differentiating the governing discrete equations and then solving a system of linear equations for the entries in the Jacobian matrix. We study the symplectic nature of the one-step algorithms for both linear and nonlinear conservative systems. It is important to note that the one-step methods developed in this paper are formulated on the state space. In order to check the condition for symplecticity, we need to define $\textbf{p}_k=m\textbf{v}_k$ to write these algorithms on phase space. Instead of writing the algorithms on phase space, we pick $m=1$ to simplify the expressions.
\par 
First, we check the condition for the simple harmonic oscillator with $L(q,\dot{q})=\frac{1}{2}\dot{q}^2 -\frac{1}{2}q^2$ where $q(t)$ is the displacement. For a fixed time step $\dt$, the Jacobians for one-step methods are 
\begin{equation*}
     \frac{\partial \phi_{\dt,V}}{\partial \textbf{y}_k} = \begin{bmatrix}
  \frac{(44\dt^5 + 143\dt^4 - 700\dt^3 + 189\dt^2 + 1176\dt - 882) }{(7\dt(\dt^4 + 9\dt^2 + 210))}&  \frac{(8\dt^5 + 33\dt^4 - 224\dt^3 - 217\dt^2 + 1568\dt + 294) }{(7(\dt^4 + 9\dt^2 + 210))} \\ \frac{-(66\dt^5 + 169\dt^4 - 434\dt^3 + 1092\dt^2 + 588\dt + 1764) }{(14\dt(\dt^4 + 9\dt^2 + 210))}&  \frac{-(12\dt^5 + 39\dt^4 - 224\dt^3 - 56\dt^2 + 784\dt - 588) }{(14\dt(\dt^4 + 9\dt^2 + 210))} \end{bmatrix},
\end{equation*}
\begin{equation*}
    \frac{\partial \phi_{\dt,G}}{\partial \textbf{y}_k} = \begin{bmatrix}
  \frac{ (2( \dt^4 - 52\dt^2 + 120))}{\dt^4 + 16\dt^2 + 240}&  \frac{(\dt(\dt^4 - 132\dt^2 + 1440)) }{6(\dt^4 + 16\dt^2 + 240)} \\ \frac{ -(\dt(\dt^4 - 130\dt^2 + 1200 )) }{5(\dt^4 + 16\dt^2 + 240)}&  \frac{-(\dt^6 - 240\dt^4 + 6240\dt^2  - 14400) }{60(\dt^4 + 16\dt^2 + 240)} \end{bmatrix}. 
\end{equation*}
For this linear  dynamical system, Jacobians from both variational and Galerkin one-step methods satisfy the condition for symplecticity. For a general mechanical system with a nonlinear potential energy $U(q)$ and Lagrangian $L=\frac{1}{2}m\dot{q}^2 -U(q)$, we find that \emph{none} of the one-step methods satisfy the required condition for symplecticity. It is important to note that the above condition is only to check whether the algorithms preserve the canonical symplectic form. In fact, at present, one can only check whether a given integration scheme exhibits a specific symplectic structure; one cannot determine whether any such symplectic structure exists, in general. In the past, some of the well-known methods such as Newmark methods \cite{kane2000variational} have been shown to preserve a noncanonical symplectic form via nonlinear transformations but it is not generally known how to test for the existence of a noncanonical symplectic form for a given algorithm. 
\section{Numerical Results}
\label{s:5}
In this section, we consider three examples to study the numerical performance of the proposed Galerkin and variational one-step methods. We first consider a nonlinear conservative system to demonstrate the good energy performance of the proposed methods. We also present an order analysis study to understand the convergence behavior of one-step variational and Galerkin methods. We then consider the Duffing oscillator to investigate the numerical behavior of the proposed one-step methods in the presence of dissipative forces. Finally, we consider a nonlinear aeroelastic system to study how the proposed methods perform for a coupled nonlinear dynamical system.
\subsection{Particle in a Double-well Potential}
\begin{figure}
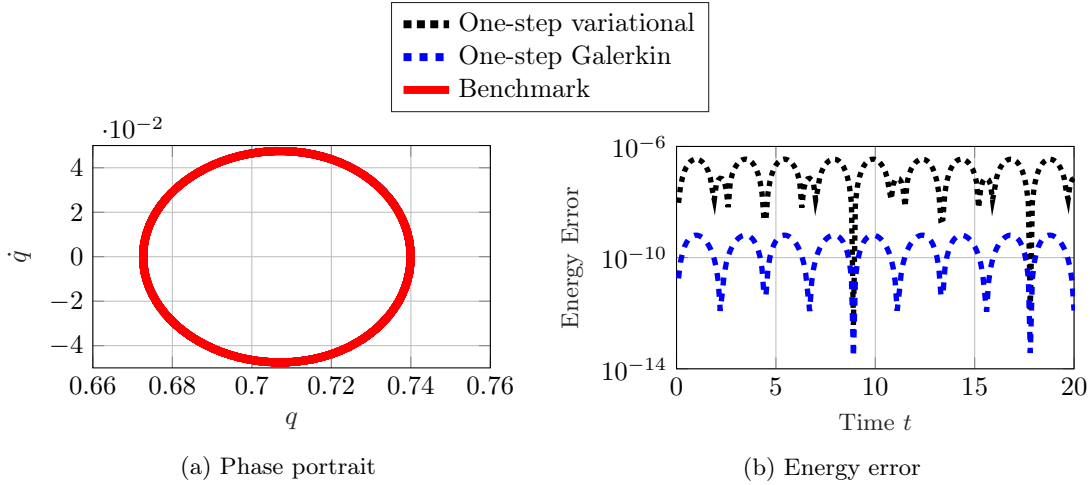

\captionsetup[subfigure]{oneside,margin={1.8cm,0 cm}}
    \begin{subfigure}{.37\textwidth}
       \setlength\fheight{5.5 cm}
        \setlength\fwidth{\textwidth}
\input{Figures/tikz/cons_pp_circ.tex}
 \caption{Phase portrait}
\label{fig:cons_pp_circ}
    \end{subfigure}
    \hspace{1.5cm}
\begin{subfigure}{.37\textwidth}
       \setlength\fheight{5.5 cm}
        \setlength\fwidth{\textwidth}
\raisebox{-58mm}{\input{Figures/tikz/cons_E_circ.tex}}
 \caption{Energy error}
\label{fig:cons_energy_circ}
    \end{subfigure}
    \caption{{\Ra{Comparison between the two one-step methods for $q(0)=0.74, \ \dot{q}(0)=0$.}}}
\label{fig:circ}
\end{figure}

\begin{figure}
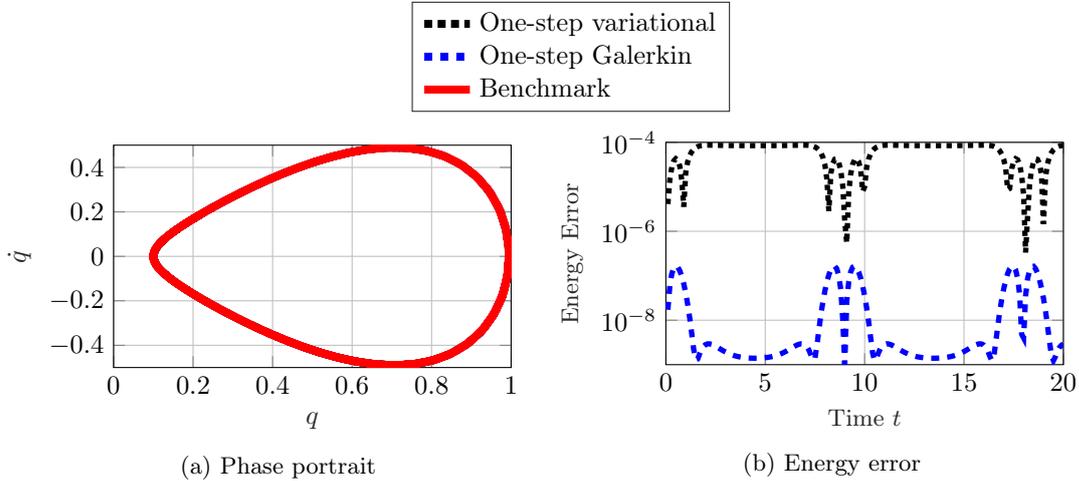

\captionsetup[subfigure]{oneside,margin={1.8cm,0 cm}}
    \begin{subfigure}{.37\textwidth}
           \setlength\fheight{5.5 cm}
           \setlength\fwidth{\textwidth}
\input{Figures/tikz/cons_pp_oval.tex}
 \caption{Phase portrait}
\label{fig:cons_pp_oval}
    \end{subfigure}
    \hspace{1.5cm}
    \begin{subfigure}{.37\textwidth}
           \setlength\fheight{5.5 cm}
           \setlength\fwidth{\textwidth}
\raisebox{-57mm}{\input{Figures/tikz/cons_E_oval.tex}}
 \caption{Energy error}
\label{fig:cons_energy_oval}
    \end{subfigure}
    \caption{{\Ra{Comparison between the two one-step methods for $q(0)=0.995, \ \dot{q}(0)=0$.}}}
\label{fig:oval}
\end{figure}

In this subsection, we apply the proposed methods to a particle in a double-well potential with Lagrangian 
\begin{equation}
    L(q,\dot{q})= \frac{1}{2}m\dot{q}^2 - \frac{1}{2} \left( q^4 -q^2 \right),
\end{equation}
The Euler-Lagrange equation for this conservative system is given by
\begin{equation}
    m\ddot{q} - q + 2q^3 = 0 .
\end{equation}
We have compared the numerical results for $m=1$ and fixed time step $\dt=0.1$ for two initial conditions. The phase portrait comparisons in Figure \ref{fig:cons_pp_circ} and Figure \ref{fig:cons_pp_oval}  show how both variational and Galerkin methods agree with the benchmark solution for both initial conditions. The corresponding energy error plots in Figure \ref{fig:cons_energy_circ} and Figure \ref{fig:cons_energy_oval} demonstrate the bounded energy error for the one-step methods. For both cases, the energy performance for the Galerkin method is substantially better than for the variational method. The energy error comparison in Figure \ref{fig:cons_energy_circ} shows that the one-step variational method has energy error magnitude around $10^{-7}$ whereas the Galerkin method has energy error around $10^{-10}$. Similarly, in Figure \ref{fig:cons_energy_oval} the variational method has energy error magnitude around $10^{-5}$ and the Galerkin method has energy error around $10^{-8}$.
\par
\begin{figure}[h]
\captionsetup[subfigure]{oneside,margin={2cm,0 cm}}
\hspace{-0.75cm}
    \begin{subfigure}{.37\textwidth}
       \setlength\fheight{5.5 cm}
        \setlength\fwidth{\textwidth}
%
%
\definecolor{mycolor1}{rgb}{0.85098,0.32549,0.09804}%
\begin{tikzpicture}

\begin{axis}[%
width=0.951\fheight,
height=0.536\fheight,
at={(0\fheight,0\fheight)},
scale only axis,
xmode=log,
xmin=0.001,
xmax=1,
xminorticks=true,
xlabel style={font=\color{white!15!black}},
xlabel={\small Time step $\Delta t$ },
ymode=log,
ymin=1e-16,
ymax=1,
ytick={1,1e-4,1e-8,1e-12,1e-16},
yminorticks=true,
ylabel style={font=\color{white!15!black}},
ylabel={\small Maximum trajectory error},
axis background/.style={fill=white},
xmajorgrids,
ymajorgrids,
legend style={at={(0.75,1.2)}, anchor=south west, legend cell align=left, align=left, draw=white!15!black,font=\small}
]
\addplot [color=red, dotted, line width=1.5pt, mark size=3.0pt, mark=diamond, mark options={solid, red}]
  table[row sep=crcr]{%
0.005	2.12816431144347e-11\\
0.01	1.49075751743055e-10\\
0.05	1.78931645056223e-08\\
0.1	2.49853970402113e-07\\
0.5	0.000224932920370513\\
};
\addlegendentry{Variational integrator}

\addplot [color=black, dotted, line width=1.5pt, mark size=3.0pt, mark=x, mark options={solid, black}]
  table[row sep=crcr]{%
0.005	5.90824276613943e-07\\
0.01	2.36326111879847e-06\\
0.05	5.9053482203697e-05\\
0.1	0.000235852294232175\\
0.5	0.00558286247047246\\
};
\addlegendentry{One-step variational}

\addplot [color=blue, dotted, line width=1.5pt, mark size=3.0pt, mark=o, mark options={solid, blue}]
  table[row sep=crcr]{%
0.005	2.65287791734181e-12\\
0.01	4.25062207654037e-11\\
0.05	2.45573754664718e-08\\
0.1	3.92752191991086e-07\\
0.5	0.000239003164719031\\
};
\addlegendentry{One-step Galerkin}

\addplot [color=purple, line width=1.5pt]
  table[row sep=crcr]{%
0.001	1e-16\\
1	0.0001\\
};
\addlegendentry{Slope 4}

\addplot [color=mycolor1, line width=1.5pt]
  table[row sep=crcr]{%
0.001	1e-06\\
1	1\\
};
\addlegendentry{Slope 2}

\end{axis}

\begin{axis}[%
width=1.227\fheight,
height=0.658\fheight,
at={(-0.16\fheight,-0.072\fheight)},
scale only axis,
xmin=0,
xmax=1,
ymin=0,
ymax=1,
axis line style={draw=none},
ticks=none,
axis x line*=bottom,
axis y line*=left
]
\end{axis}
\end{tikzpicture}%
 \caption{Trajectory error}
\label{fig:conv_q}
    \end{subfigure}
   \hspace{1.5cm}
    \begin{subfigure}{.37\textwidth}
           \setlength\fheight{5.5 cm}
           \setlength\fwidth{\textwidth}
\raisebox{-68mm}{
%
%
\definecolor{mycolor1}{rgb}{0.85098,0.32549,0.09804}%
\begin{tikzpicture}

\begin{axis}[%
width=0.951\fheight,
height=0.536\fheight,
at={(0\fheight,0\fheight)},
scale only axis,
xmode=log,
xmin=0.001,
xmax=1,
xminorticks=true,
xlabel style={font=\color{white!15!black}},
xlabel={\small Time step $\Delta t$ },
ymode=log,
ymin=1e-16,
ymax=1,
yminorticks=true,
ytick={1,1e-4,1e-8,1e-12,1e-16},
ylabel style={font=\color{white!15!black}},
ylabel={\small Maximum velocity error },
axis background/.style={fill=white},
xmajorgrids,
ymajorgrids,
legend style={at={(0.556,0.178)}, anchor=south west, legend cell align=left, align=left, draw=white!15!black,font=\tiny}
]
\addplot [color=red, dotted, line width=1.5pt, mark size=3.0pt, mark=diamond, mark options={solid, red}]
  table[row sep=crcr]{%
0.005	3.09612249033142e-11\\
0.01	2.2351859035527e-10\\
0.05	2.77475888489989e-08\\
0.1	3.1679820318084e-07\\
0.5	0.000287253162452705\\
};

\addplot [color=black, dotted, line width=1.5pt, mark size=3.0pt, mark=x, mark options={solid, black}]
  table[row sep=crcr]{%
0.005	7.853363391307e-07\\
0.01	3.14128730408408e-06\\
0.05	7.84766281019421e-05\\
0.1	0.000312879754914558\\
0.5	0.00720053405269558\\
};

\addplot [color=blue, dotted, line width=1.5pt, mark size=3.0pt, mark=o, mark options={solid, blue}]
  table[row sep=crcr]{%
0.005	3.5086950705976e-12\\
0.01	5.62867304848558e-11\\
0.05	3.28213047430351e-08\\
0.1	5.23594693768047e-07\\
0.5	0.000313763432839876\\
};

\addplot [color=purple, line width=1.5pt]
  table[row sep=crcr]{%
0.001	1e-16\\
1	0.0001\\
};

\addplot [color=mycolor1, line width=1.5pt]
  table[row sep=crcr]{%
0.001	1e-06\\
1	1\\
};

\end{axis}

\begin{axis}[%
width=1.227\fheight,
height=0.658\fheight,
at={(-0.16\fheight,-0.072\fheight)},
scale only axis,
xmin=0,
xmax=1,
ymin=0,
ymax=1,
axis line style={draw=none},
ticks=none,
axis x line*=bottom,
axis y line*=left
]
\end{axis}
\end{tikzpicture}%
}
 \caption{Velocity error}
\label{fig:conv_v}
    \end{subfigure} \\
    \centering
 \begin{subfigure}{.37\textwidth}
           \setlength\fheight{5.5 cm}
           \setlength\fwidth{\textwidth}
%
%
\definecolor{mycolor1}{rgb}{0.85098,0.32549,0.09804}%
\begin{tikzpicture}

\begin{axis}[%
width=0.951\fheight,
height=0.536\fheight,
at={(0\fheight,0\fheight)},
scale only axis,
xmode=log,
xmin=0.001,
xmax=1,
xminorticks=true,
xlabel style={font=\color{white!15!black}},
xlabel={\small Time step $\Delta t$},
ymode=log,
ymin=1e-16,
ymax=1,
ytick={1,1e-4,1e-8,1e-12,1e-16},
yminorticks=true,
ylabel style={font=\color{white!15!black}},
ylabel={\small Maximum energy error },
axis background/.style={fill=white},
xmajorgrids,
ymajorgrids,
legend style={at={(1.05,0.214)}, anchor=south west, legend cell align=left, align=left, draw=white!15!black,font=\tiny}
]
\addplot [color=red, dotted, line width=1.5pt, mark size=3.0pt, mark=diamond, mark options={solid, red}]
  table[row sep=crcr]{%
0.001	2.6703479760748e-08\\
0.005	6.67370982254668e-07\\
0.01	2.66890448880064e-06\\
0.05	6.63746192823278e-05\\
0.1	0.000261175102381939\\
0.5	0.00366937644697657\\
};

\addplot [color=black, dotted, line width=1.5pt, mark size=3.0pt, mark=x, mark options={solid, black}]
  table[row sep=crcr]{%
0.001	8.03622248260471e-09\\
0.005	2.00843990917632e-07\\
0.01	8.03502005545381e-07\\
0.05	2.00717790069145e-05\\
0.1	8.0130979268775e-05\\
0.5	0.00187160880822741\\
};

\addplot [color=blue, dotted, line width=1.5pt, mark size=3.0pt, mark=o, mark options={solid, blue}]
  table[row sep=crcr]{%
0.001	2.37049962992231e-14\\
0.005	1.53518725509727e-11\\
0.01	5.47705604828086e-11\\
0.05	9.5363903995782e-09\\
0.1	1.53726777252228e-07\\
0.5	0.000129804906132569\\
};

\addplot [color=purple, line width=1.5pt]
  table[row sep=crcr]{%
0.001	1e-15\\
1	0.001\\
};

\addplot [color=mycolor1, line width=1.5pt]
  table[row sep=crcr]{%
0.001	1e-07\\
1	0.1\\
};

\end{axis}

\begin{axis}[%
width=1.227\fheight,
height=0.658\fheight,
at={(-0.16\fheight,-0.072\fheight)},
scale only axis,
xmin=0,
xmax=1,
ymin=0,
ymax=1,
axis line style={draw=none},
ticks=none,
axis x line*=bottom,
axis y line*=left
]
\end{axis}
\end{tikzpicture}%
 \caption{Energy error}
\label{fig:conv_E}
    \end{subfigure}
    \caption{{\Ra{Convergence analysis of maximum temporal error in configuration, velocity and energy for the nonlinear conservative system.}}}
\label{fig:cons_conv}
\end{figure}
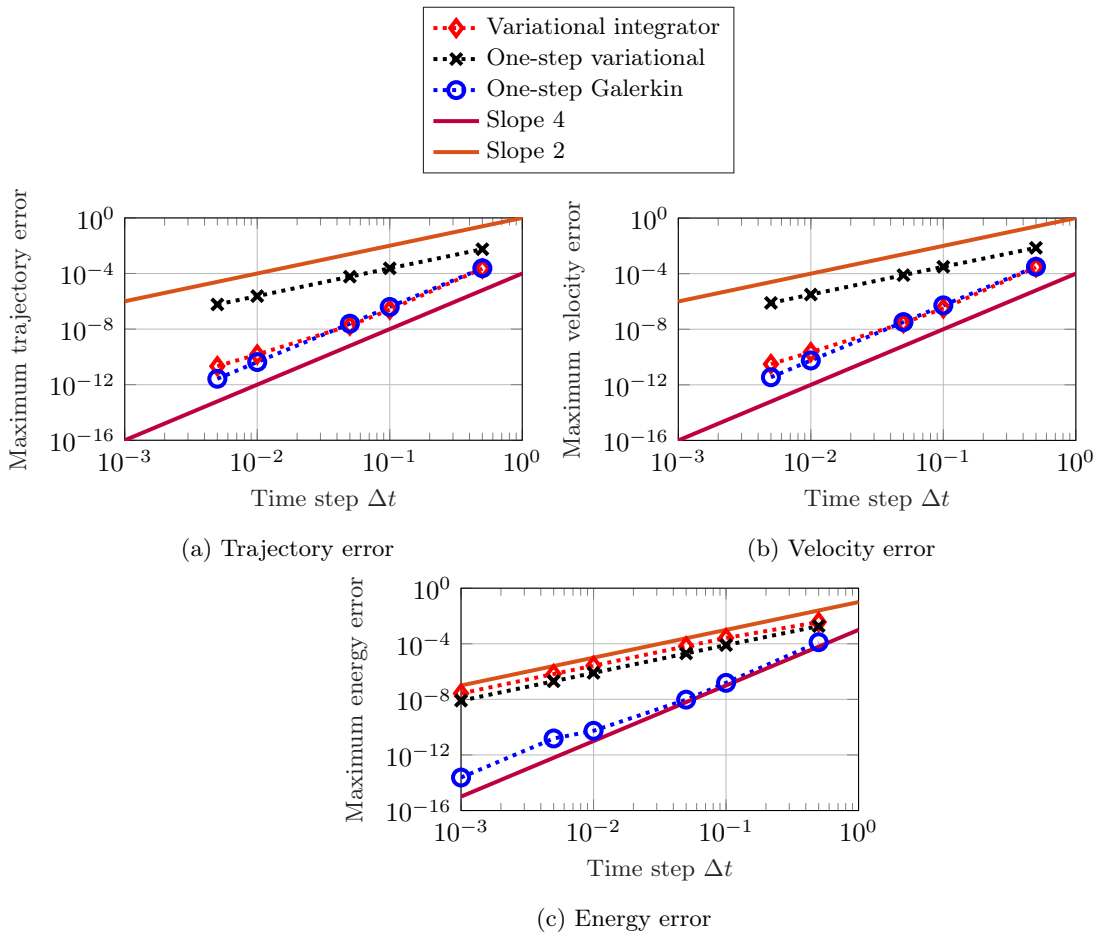

In Figure \ref{fig:cons_conv}, we have studied the numerical behavior of these algorithms for different fixed time step values to understand how the configuration, velocity and energy error values decrease with decrease in the step size. We have also performed the convergence analysis for variational integrators derived from the discrete mechanics framework to understand how they compare to the proposed one-step methods. Since both one-step methods proposed in this work require the solution of two coupled nonlinear equations at each time step, we have implemented variational integrators that require solving two nonlinear equations. Thus, we have considered a quadratic trajectory over one time step in the discrete mechanics framework by introducing an interior point \cite{marsden2001discrete} and the resulting variational integrator leads to two coupled nonlinear implicit equations at each time step. 
\par
The configuration and velocity convergence plots in Figure \ref{fig:conv_q} and Figure \ref{fig:conv_v} show that the one-step variational method has second order convergence whereas both the variational integrator and the one-step Galerkin method have fourth order convergence. For the energy error, the one-step variational and the variational integrator show second order convergence, whereas the one-step Galerkin method shows fourth order convergence. Thus, the one-step Galerkin approach gives better trajectory and energy performance than both the one-step variational method and the variational integrator.
\subsection{Duffing Oscillator}
\begin{figure}
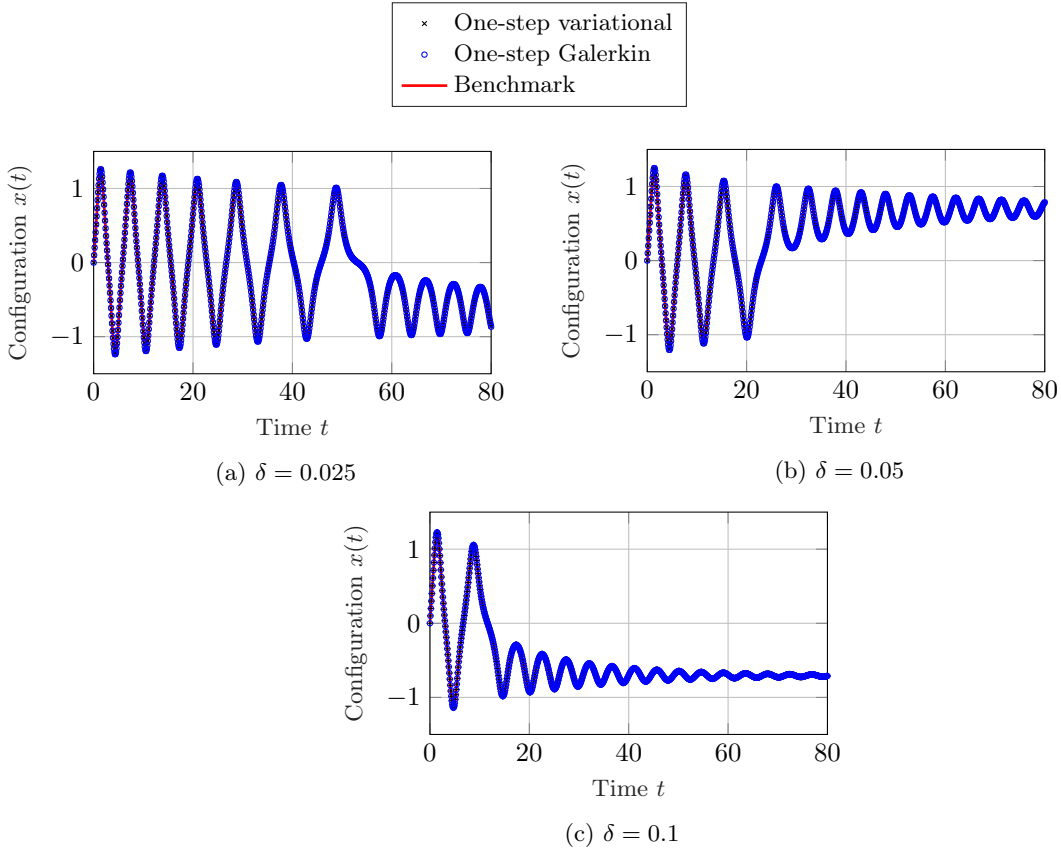

\captionsetup[subfigure]{oneside,margin={2cm,0 cm}}
\hspace{-0.75cm}
    \begin{subfigure}{.37\textwidth}
       \setlength\fheight{5.5 cm}
        \setlength\fwidth{\textwidth}
\input{Figures/tikz/cd_0025_comp.tex}
 \caption{$\delta=0.025$}
\label{fig:cd_1}
    \end{subfigure}
   \hspace{1.5cm}
    \begin{subfigure}{.37\textwidth}
           \setlength\fheight{5.5 cm}
           \setlength\fwidth{\textwidth}
\raisebox{-58mm}{\input{Figures/tikz/cd_005_comp.tex}}
 \caption{$\delta=0.05$}
\label{fig:cd_2}
    \end{subfigure} \\
    \centering
 \begin{subfigure}{.37\textwidth}
           \setlength\fheight{5.5 cm}
           \setlength\fwidth{\textwidth}
\input{Figures/tikz/cd_01_comp.tex}
 \caption{$\delta=0.1$}
\label{fig:cd_3}
    \end{subfigure}
    \caption{{\Ra{Duffing oscillator numerical simulation with $\alpha=-1,\beta=2$, and time step $\dt=0.1$.}}}\label{fig:cd}
\end{figure}

In this subsection, we study the numerical performance of these algorithms in the presence of dissipation. The governing second-order, nonlinear differential equation for the Duffing oscillator is 
\begin{equation}
 \ddot{x} + \delta \dot{x} + \alpha x + \beta x^3 = 0,
\end{equation}
where $x(t)$ is the displacement at time $t$, $\delta$ is the linear damping, $\alpha$ is the linear stiffness, and $\beta$ is the nonlinear stiffness coefficient. The Lagrangian and external forcing for this system are 
\begin{equation}
    L(x,\dot{x})=\frac{1}{2}\dot{x}^2  - \frac{1}{2}\alpha x^2 -\frac{1}{4}\beta x^4, \quad \quad f(\dot{x},t) = -\delta \dot{x} .
\end{equation}
We have fixed the stiffness parameters to $\alpha=1$ and $\beta=0.5$ and studied this dissipative nonlinear dynamical system for three cases, i.e. $\delta \in \{0.025,0.05,0.1\}$.  With these specific parameter values, the Duffing oscillator can be thought of as the double-well potential system with dissipation. The numerically computed trajectories from the one-step methods are compared with the benchmark solution in Figure \ref{fig:cd}. The plots in Figure \ref{fig:cd} demonstrate that both methods are able to capture the dissipation effect accurately. In fact, the discrete trajectories are indistinguishable from the benchmark solution for all three cases.
\par
\begin{figure}[h]
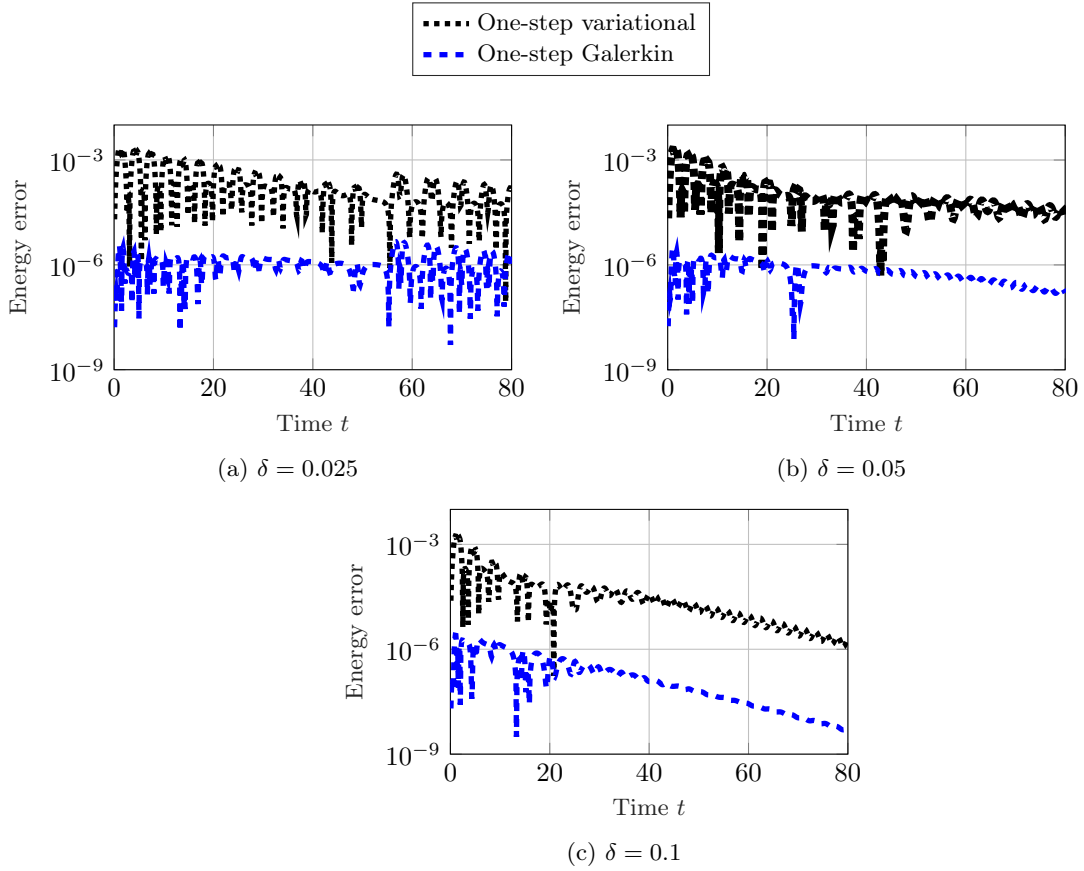

\captionsetup[subfigure]{oneside,margin={2cm,0 cm}}
\hspace{-0.75cm}
    \begin{subfigure}{.37\textwidth}
       \setlength\fheight{5.5 cm}
        \setlength\fwidth{\textwidth}
\input{Figures/tikz/cd_0025_Ecomp.tex}
 \caption{$\delta=0.025$}
\label{fig:cd_E1}
    \end{subfigure}
   \hspace{1.5cm}
    \begin{subfigure}{.37\textwidth}
           \setlength\fheight{5.5 cm}
           \setlength\fwidth{\textwidth}
\raisebox{-58mm}{\input{Figures/tikz/cd_005_Ecomp.tex}}
 \caption{$\delta=0.05$}
\label{fig:cd_E2}
    \end{subfigure} \\
    \centering
 \begin{subfigure}{.37\textwidth}
           \setlength\fheight{5.5 cm}
           \setlength\fwidth{\textwidth}
\input{Figures/tikz/cd_01_Ecomp.tex}
 \caption{$\delta=0.1$}
\label{fig:cd_E3}
    \end{subfigure}
    \caption{{\Ra{Energy error behavior for Duffing oscillator numerical simulation with $\alpha=-1,\beta=2$, and time step $\dt=0.1$.}}}\label{fig:cd_E}
\end{figure}

The energy error plots in Figure \ref{fig:cd_E} compare the energy performance for the one-step methods and the Galerkin approach outperforms the variational approach in all three cases. For all three cases, the energy error for the one-step variational method starts around $10^{-4}$ whereas the one-step Galerkin methods exhibit energy error around $10^{-6}$. The energy error for both methods decreases over the time due to the presence of dissipative forces and the rate of decrease in energy error increases with increasing values of the damping parameter $\delta$. This decrease in energy error is seen clearly in Figure \ref{fig:cd_E2} and Figure \ref{fig:cd_E3} for higher $\delta$ values.
\subsection{Aeroelastic System}
\begin{figure}[h]
\centering
 \includegraphics[width=0.6\textwidth]{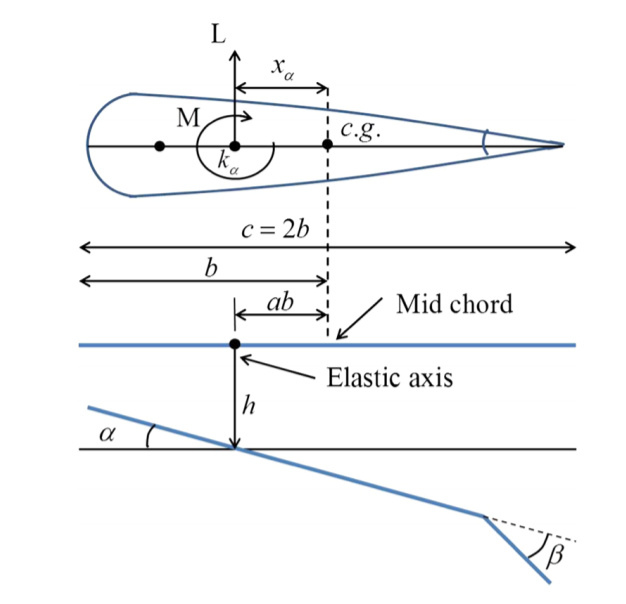}
 \caption{Sketch of a two degrees of freedom aeroelastic section model.}
 \label{fig: aeroelas}
\end{figure}

In this subsection, we consider the open-loop behavior of the nonlinear aeroelastic system studied by Shukla and Patil \cite{shukla2017nonlinear}. As shown in Figure \ref{fig: aeroelas}, the model contains a flat plate supported by a linear spring in the plunge degree of freedom and the cubic nonlinear spring in the pitch degree of freedom. The flat plate is free to move up and down along the plunge degree of freedom and to rotate about the pitch degree of freedom. The Lagrangian for this system is given by
\begin{equation}
    L(h, \dot{h},\alpha,\dot{\alpha})= m_{\rm T} \dot{h}^2 + I_{\alpha}\dot{\alpha}^2 + m_{\rm W}x_{\alpha}\dot{h}\dot{\alpha} - \frac{1}{2}k_hh^2 -\frac{1}{2}k_{\alpha_0}\alpha^2 -  \frac{1}{3}k_{\alpha_1}\alpha^3 - \frac{1}{4}k_{\alpha_2}\alpha^4  ,
\end{equation}
where $m_{\rm W}$ is the mass of the wing and $m_{\rm T}$ is the total mass of the aeroelastic system. The parameter $I_{\alpha}$ represents the moment of inertia about the elastic axis. The terms $k_h$ and $k_{\alpha_{\{0,1,2\}}}$ are the stiffness functions along the plunge and pitch degrees of freedom respectively.
The external nonconservative forces are 
\begin{equation}
    f_{h}= - c_h \dot{h} + \rho U^2bC_{L_{\alpha}} \alpha_{\rm eff},
\end{equation}
\begin{equation}
    f_{\alpha} = -c_{\alpha}\dot{\alpha} + \rho U^2b^2C_{M_{\alpha}} \alpha_{\rm eff},
\end{equation}
where $c_h$ and $c_{\alpha}$ are damping coefficients, $C_{L_{\alpha}}$ and $C_{M_{\alpha}}$ are the derivatives of the lift and moment coefficients, and  $\alpha_{\rm eff}=\left( \alpha + \frac{\dot{h}}{U} + \left(\frac{1}{2}-a \right)b \frac{\dot{\alpha}}{U} \right)$ is the effective angle of attack. The equations of motion for this aeroelastic system are 
\begin{equation}
    m_{\rm T}\ddot{h} + m_{\rm W}x_{\alpha}b\ddot{\alpha} + c_h \dot{h} + k_hh -  \rho U^2bC_{L_{\alpha}} \alpha_{\rm eff} = 0,
\end{equation}
\begin{equation}
    I_{\alpha}\ddot{\alpha} + m_{\rm W}x_{\alpha}b\ddot{h} + c_{\alpha} \dot{\alpha} + k_{\alpha}(\alpha)\alpha + \rho U^2b^2C_{M_{\alpha}} \alpha_{\rm eff}   = 0,
\end{equation}
where $k(\alpha)=k_{\alpha_0} + k_{\alpha_1}\alpha + k_{\alpha_2}\alpha^2$.
\par
\begin{figure}
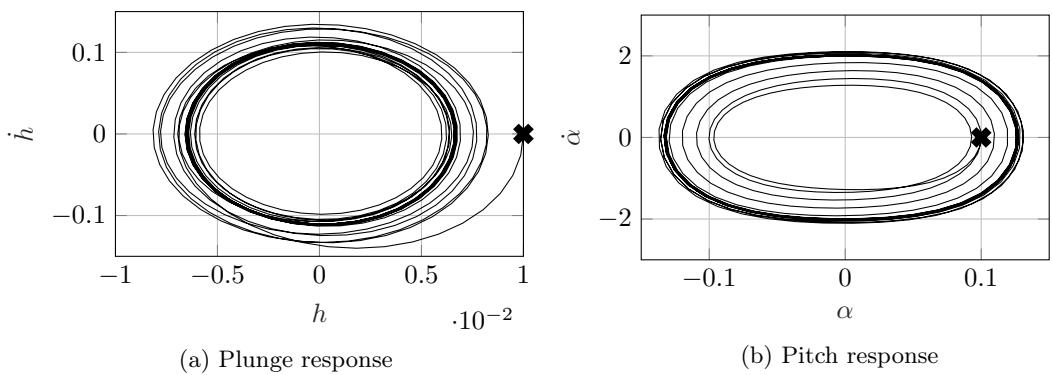

\captionsetup[subfigure]{oneside,margin={2cm,0 cm}}
    \begin{subfigure}{.37\textwidth}
       \setlength\fheight{5.5 cm}
        \setlength\fwidth{\textwidth}
\input{Figures/tikz/var_h_pp.tex}
 \caption{Plunge response}
\label{fig: var_h_pp}
    \end{subfigure}
   \hspace{1.5cm}
    \begin{subfigure}{.37\textwidth}
           \setlength\fheight{5.5 cm}
           \setlength\fwidth{\textwidth}
\input{Figures/tikz/var_a_pp.tex}
 \caption{Pitch response}
\label{fig: var_a_pp}
   \end{subfigure}
    \caption{{\Ra{Subcritical LCO simulation using the one-step variational method.}}}\label{fig:pp_aero_var}
\end{figure}
\begin{figure}
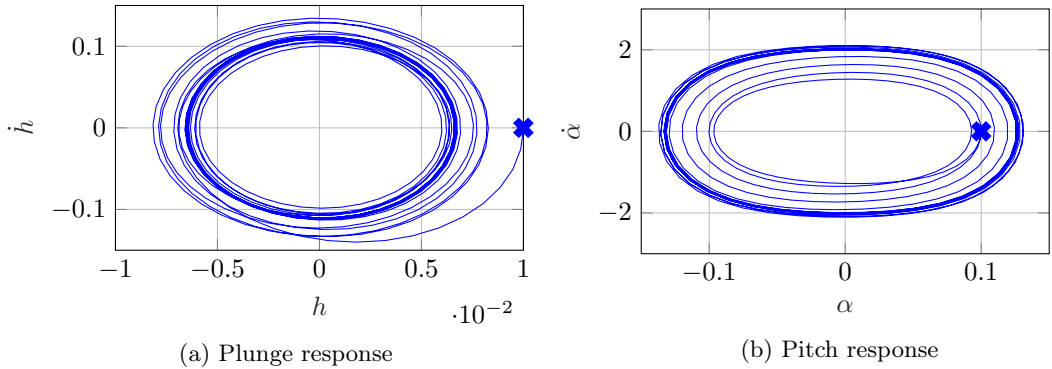

\captionsetup[subfigure]{oneside,margin={2cm,0 cm}}
    \begin{subfigure}{.37\textwidth}
       \setlength\fheight{5.5 cm}
        \setlength\fwidth{\textwidth}
\input{Figures/tikz/gal_h_pp.tex}
 \caption{Plunge response}
\label{fig: gal_h_pp}
    \end{subfigure}
   \hspace{1.5cm}
    \begin{subfigure}{.37\textwidth}
           \setlength\fheight{5.5 cm}
           \setlength\fwidth{\textwidth}
\input{Figures/tikz/gal_a_pp.tex}
 \caption{Pitch response}
\label{fig: gal_a_pp}
   \end{subfigure}
    \caption{{\Ra{Subcritical LCO simulation using the one-step Galerkin method.}}}\label{fig:pp_aero_gal}
\end{figure}

We have studied the dynamic behavior of the nonlinear aeroelastic system for initial conditions $(h_0,\dot{h}_0,\alpha_0,\dot{\alpha}_0)=(0.01,0,0.1,0)$ at freestream velocity $U=0.9U_f$ where $U_f$ is the linear flutter velocity. The phase space plots given in Figure \ref{fig:pp_aero_var} and Figure \ref{fig:pp_aero_gal} clearly demonstrate how both one-step methods capture the subcritical limit cycle oscillations (LCOs) accurately.
\par
\begin{figure}[h]
\centering
    \setlength\fheight{10.5 cm}
        \setlength\fwidth{\textwidth}
\input{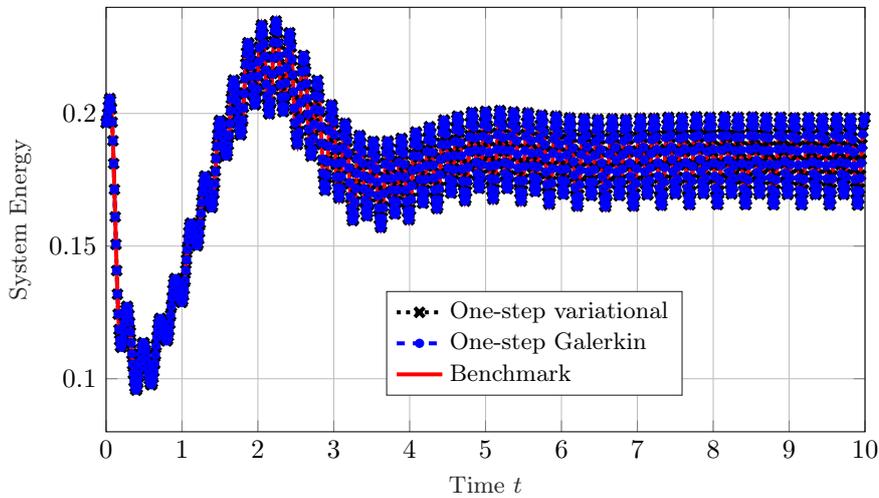}
 \caption{{\Ra{Total energy of the aeroelastic system.}}}
 \label{fig:ener_aero}
\end{figure}
\begin{figure}[h]
\centering
    \setlength\fheight{10.5 cm}
        \setlength\fwidth{\textwidth}
\input{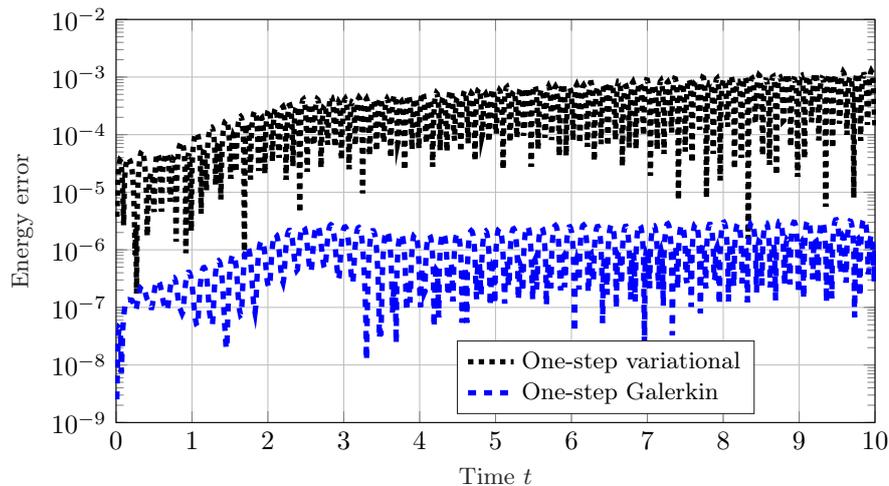}
  \caption{{\Ra{Energy error comparison between one-step Galerkin and variational methods.}}}
    \label{fig:err}
\end{figure}

The energy plot in Figure \ref{fig:ener_aero} shows how the total energy of the aeroelastic system evolves over time. Initially there is a sharp decrease in energy followed by an increase, and eventually when the system exhibits periodic motion with constant amplitude the total energy oscillates around a fixed value. As shown in Figure \ref{fig:ener_aero}, both one-step methods track the change in energy accurately for the nonlinear aeroelastic system. The energy error comparison in Figure \ref{fig:err} demonstrates how the Galerkin approach has better energy behavior than the one-step variational approach. The one-step variational method has energy error magnitude around $10^{-4}$ whereas the Galerkin method has energy error around $10^{-6}$.
\section{Conclusions}
\label{s:6}
In this paper we have developed Hermite polynomial based one-step variational and Galerkin methods for mechanical systems with external forcing. We have utilized cubic Hermite polynomials over one time step for discretization and the resulting numerical algorithms are continuous in both configuration and velocity. We also demonstrated an approach to obtain one-step methods using higher-order Hermite polynomials. We showed that both one-step methods are symplectic for linear dynamical systems but they do not preserve the canonical symplectic form for general nonlinear dynamical systems. We also investigated the linear stability of the proposed one-step methods and both methods exhibit excellent stability for large time steps.
\par 
We have studied the numerical behavior of these algorithms through three different numerical examples. The energy performance and convergence analysis results for the conservative example showed how both one-step methods achieve good numerical performance by obtaining $\mathcal{C}^1$--continuous trajectories. We have also presented results for a dissipative system and the numerical plots show that both one-step methods capture the effect of the dissipative forces accurately over long time intervals. Finally, numerical studies for the coupled aeroelastic system show how both one-step methods capture the limit cycle oscillations accurately. The numerical results from all three examples showed that the Galerkin approach has significantly better energy behavior. In fact, the one-step Galerkin method is also better than the variational integrator, with same nonlinear equations per time step, in terms of the energy accuracy.  
\par 
{\Rb{Future research directions motivated by this work are: obtaining theoretical results about the convergence behavior and geometric properties of the proposed one-step methods; investigating the connection between the one-step Galerkin methods and energy-momentum integrators; and applying these one-step methods to discretizations of infinite-dimensional systems.}}
\\ 
\\
\textbf{Declaration of interest:} None.
\\ 
\\
\\
\\
\textbf{Funding:} This material is based upon work supported by the National Science Foundation under Grant No. 1826152.
\bibliographystyle{vancouver}
\bibliography{main}

\end{document}